\documentclass[a4paper,reqno,10pt,oneside]{amsart}
\usepackage{amsmath}
\usepackage[left=2.61cm,right=2.61cm,top=2.72cm,bottom=2.72cm]{geometry}
\usepackage[colorlinks=true,urlcolor=blue,
citecolor=red,linkcolor=blue,linktocpage,pdfpagelabels,
bookmarksnumbered,bookmarksopen]{hyperref}
\usepackage[english]{babel}
\usepackage{graphicx}
\usepackage[toc,page,header]{appendix}
\usepackage{mathrsfs}
\usepackage{amssymb,amsmath, amsfonts, amsthm}
\usepackage{latexsym}
\usepackage{enumitem}

\allowdisplaybreaks[1]

\numberwithin{equation}{section}

\DeclareMathOperator{\divergence}{div}
\DeclareMathOperator{\diver}{div}
\DeclareMathOperator{\Ric}{Ric}

\DeclareMathOperator{\supp}{supp}

\renewcommand{\(}{\left(}
\renewcommand{\)}{\right)}
\renewcommand{\[}{\left[}
\renewcommand{\]}{\right]}

\newtheorem{theorem}{Theorem}[section]
\newtheorem{proposition}[theorem]{Proposition}
\newtheorem{corollary}[theorem]{Corollary}
\newtheorem{lemma}[theorem]{Lemma}
\theoremstyle{definition}
\newtheorem{remark}[theorem]{Remark}
\theoremstyle{definition}

\theoremstyle{definition}

\newcommand{\vp}{\varphi}

\renewcommand{\le}{\leqslant}
\renewcommand{\ge}{\geqslant}

\newcommand{\beq}{\begin{equation}}
\newcommand{\eeq}{\end{equation}}
\newcommand{\beqs}{\begin{equation*}}
\newcommand{\eeqs}{\end{equation*}}
\newcommand{\beqn}{\begin{eqnarray}}
\newcommand{\eeqn}{\end{eqnarray}}
\newcommand{\beqns}{\begin{eqnarray*}}
\newcommand{\eeqns}{\end{eqnarray*}}
\newcommand{\bdoc}{\begin{document}}
\newcommand{\edoc}{\end{document}}
\newcommand{\be}{\begin{enumerate}}
\newcommand{\ee}{\end{enumerate}}
\newcommand{\bdescr}{\begin{description}}
\newcommand{\edescr}{\end{description}}
\newcommand{\ba}{\begin{array}}
\newcommand{\ea}{\end{array}}
\newcommand{\intR}{\int_{\mathbb R^N}}
\newcommand{\R}{\mathbb R}
\newcommand{\RN}{\mathbb{R}^N}
\newcommand{\B}{\mathbb B}
\renewcommand{\H}{\mathcal H}
\renewcommand{\L}{\mathbb L}
\newcommand{\parallelsum}{\mathbin{\!/\mkern-5mu/\!}}
\newcommand{\e}{\varepsilon}
\newcommand{\SD}{\Sigma_D}
 \renewcommand{\(}{\left(}
\renewcommand{\)}{\right)}
\renewcommand{\[}{\left[}
\renewcommand{\]}{\right]}
\renewcommand{\appendixpagename}{\centering Appendix}

\newcommand{\todo}[1]{\text{\colorbox{yellow}{#1}}}

\begin{document}
\title[Semilinear PDEs on Riemannian manifolds :  a P-function approach]{Classification results, rigidity theorems \\ and semilinear PDEs on Riemannian manifolds : \\ a P-function approach}
%P-function for semilinear PDEs with nonnegative Ricci curvature]{Classification results via P-function for semilinear PDEs in manifolds}

\author{Giulio Ciraolo}
\address{G. Ciraolo. Dipartimento di Matematica "Federigo Enriques",
Universit\`a degli Studi di Milano, Via Cesare Saldini 50, 20133 Milano, Italy}
\email{giulio.ciraolo@unimi.it}

\author{Alberto Farina}
\address{A. Farina. LAMFA, UMR CNRS 7352, Universit\'e Picardie Jules Verne 33, rue St Leu, 
80039 Amiens, France}
\email{alberto.farina@u-picardie.fr}

\author{Camilla Chiara Polvara}
\address{C.C. Polvara. Dipartimento di Matematica "Federigo Enriques",
Universit\`a degli Studi di Milano, Via Cesare Saldini 50, 20133 Milano, Italy}
\email{camilla.polvara@unimi.it}



%\thanks{\emph{Acknowledgements.} Research partially supported by Gruppo Nazionale per l'Analisi Matematica, la Pro\-ba\-bi\-li\-t\`a e le loro Applicazioni (GNAMPA) of the Istituto Nazionale di Alta Matematica (INdAM)}

\begin{abstract}
We consider solutions to some semilinear elliptic equations on complete noncompact Riemannian manifolds and study their classification as well as the effect of their presence on the underlying manifold.
When the Ricci curvature is non-negative, we prove both the classification of positive solutions to the critical equation and the rigidity for the ambient manifold. The same results are obtained when we consider solutions to the Liouville equation on Riemannian surfaces. The results are obtained via a suitable P-function whose constancy implies the classification of both the solutions and the underlying manifold.
The analysis carried out on the P-function also makes it possible to classify non-negative solutions for subcritical equations on manifolds enjoying a Sobolev inequality and satisfying an integrability condition on the negative part of the Ricci curvature.
Some of our results are new even in the Euclidean case.
\end{abstract}

\maketitle

\section{Introduction}

In this paper we consider solutions to critical and subcritical semilinear elliptic equations on complete noncompact Riemannian manifolds. We are interested in their classification as well as on the effect of their presence on the ambient manifold.

%In this paper we consider some classical problems related to the classification of entire solutions to elliptic partial differential equations. 
%One of these problems is the classification of positive solutions of the critical Laplace equation

The prototype of one of these problems is the classification of positive solutions of the critical Laplace equation 
\begin{equation} \label{eq_crit_Eucl}
- \Delta u = u^{\frac{d+2}{d-2}} \quad \textmd{ in } \R^d \,,
\end{equation} 
with $d \geq 3$, which arises as the Euler-Lagrange equation of the $L^2$ Sobolev inequality as well as from the study of Yamabe problem, where $u^{\frac{4}{d-2}}g_{\R^d}$ provides a metric of positive and constant scalar curvature. It is well known that positive solutions to \eqref{eq_crit_Eucl} can be completely classified, see the classical papers \cite{Aubin,GNN,Obata,Talenti} where the classification was obtained under further decay  or variational assumptions on the solution. These assumptions were removed in the celebrated paper \cite{CGS} (and later in \cite{CL1991} and \cite{LiZhang}), where the classification result is obtained by only requiring the positivity of the solution. The results of \cite{CGS,CL1991,LiZhang}, whose proofs are crucially based on the moving plane method and on the Kelvin transform, tell us that positive solutions to \eqref{eq_crit_Eucl} assume the form 
\begin{equation}
u(x) = \left( \frac{1}{a+b|x-x_0|^2} \right)^{\frac{d-2}{2}} \,,
\end{equation}
for some $a,b>0$ with $1=d(d-2)ab$ and $x_0 \in \R^d$.  

A related result was obtain in \cite{GS}, where it was shown the nonexistence of positive solutions to subcritical equations 
\begin{equation} \label{eq_subcr_Rd}
- \Delta u = u^p \quad \textmd{ in } \R^d \,,
\end{equation}
where $1 \leq p < \frac{d+2}{d-2}$, with $d \geq 3$. 

When $d=2$, a different line of investigation comes from the classification of entire solutions of the Liouville equation 
\begin{equation} \label{Liouv_Rd}
- \Delta u = e^{u} \quad \textmd{ in } \R^2 \,,
\end{equation}
which also arises from the study of conformal metrics $e^u g_{\R^2}$. This problem goes back to the work of J. Liouville \cite{Liouville} who found that solutions to \eqref{Liouv_Rd} can be represented by some meromorphic functions on the complex plane $\mathbb{C}$. Several years later, in \cite{CL1991}, solutions to \eqref{Liouv_Rd} were completely classified under the finite mass condition 
\begin{equation} \label{finite-mass-Liouv_Rd}
\int_{\mathbb{R}^2} e^u\,dx<+\infty \,.
\end{equation}
Under the latter condition, the authors of \cite{CL1991} prove that solutions of \eqref{Liouv_Rd} are given by  
\begin{equation}\label{u_radial_liouv}
u(x) = \log \frac{1}{(a+b|x-x_0|^2)^2} \,,
\end{equation}
for some $a,b>0$ with $1=8ab$ and $ x_0 \in \R^2$.  The proof in \cite{CL1991} relies on the moving plane procedure and on the knowledge of the precise asymptotic behaviour of the solution at infinity (the latter being obtained thanks to the finite mass condition \eqref{finite-mass-Liouv_Rd} and a result of  \cite{BM}).

\medskip

In what follow, we first consider problems \eqref{eq_crit_Eucl} and \eqref{Liouv_Rd} on complete noncompact Riemannian manifolds with nonnegative Ricci curvature and we aim at providing a unified approach that allows us to classify both the solutions to those equations and the ambient manifold. Such results will be obtained by using a suitable $P$-function which captures all the ``geometry'' encoded in equations \eqref{eq_crit_Eucl} and \eqref{Liouv_Rd} and, whose constancy implies both the classification of the solutions and the rigidity of the underlying manifold.  
Then, by exploiting the analysis performed on the P-function, we classify non-negative solutions for subcritical equations on manifolds enjoying a Sobolev inequality and satisfying an integrability condition on the negative part of the Ricci curvature.

\medskip

Our first main result is the following.

\begin{theorem} \label{thm_critical}
Let $(M,g)$ be a complete, connected, non-compact, boundaryless Riemannian manifold 
of dimension $d \geq 3 $ and with nonnegative Ricci curvature. Let $u \in C^3(M)$ be a positive solution to 
\begin{equation} \label{eq_crit}
- \Delta u = u^{\frac{d+2}{d-2}} \qquad \textmd{ in } \quad M.
\end{equation}
Assume that one (and only one) of the following three conditions holds
\begin{itemize}
\item[(i)] $d=3,4,5$;
\item[(ii)]$d \geq 6$ and $u(x) = O(r^{-\alpha})$ with 
\begin{equation} \label{alpha_opt}
\alpha > \frac{(d-6)(d-2)}{4(d-4)} \,,
\end{equation}
where $r=r(x)$ denotes the geodesic distance from a fixed reference point of $M$;
\item[(iii)] $d\geq 3$ and $u \in L^{\frac{2d}{d-2}}(M)$.
\end{itemize}
Then $(M,g)$  is isometric to $\R^d$ with the Euclidean metric and $u$ is given by
\begin{equation}\label{u_radial_crit}
u(x) = \left( \frac{1}{a+b|x-x_0|^2} \right)^{\frac{d-2}{2}} \,,
\end{equation}
for some $a,b>0$ with $1=d(d-2)ab$ and $x_0 \in \R^d$. 
\end{theorem} 

Note that, in particular, Theorem \ref{thm_critical} tells us that the non-compact Yamabe problem on a complete, connected, boundaryless  Ricci-flat manifold of dimension $d=3,4,5$ is solvable if and only if the manifold is isometric to $\R^d$ with the Euclidean metric. 

Theorem \ref{thm_critical} recovers and improves the results available in literature in several directions. 
%Besides the classical results in the Euclidean case, our results improve those contained in the very recent papers 
%\cite{FMM} and \cite{CatinoMonticelli}.  

In \cite{FMM}, the conclusion of Theorem \ref{thm_critical} was obtained under the decay assumption $u=O(r^{-\frac{d-2}{2}}) $ at infinity (together with a further assumption either on the volume growth of $M$ or on the finiteness of the energy of $u$). Since the asymptotic assumption made in \cite{FMM} immediately implies that $u$ satisfies both (ii) and (iii), we fully recover the results in \cite{FMM}.  In \cite{CatinoMonticelli}, the case $(i)$ with $d=3$ was already obtained, while the classification result was proved for dimension $d \geq 4$ either under the finite energy condition (i.e., $u \in L^{\frac{2d}{d-2}}(M)$ and $\nabla u \in L^{2}(M)$)\footnote{The authors of \cite{CatinoMonticelli} informed us that  
item (iii) of our Theorem \ref{thm_critical} can also be deduced by combining  Lemma 2.9 and Theorem 1.2 of \cite{CatinoMonticelli}.} or under a decay assumption of the solution at infinity (more restrictive than our condition $(ii)$ when $d \geq 7$, and coinciding with $(ii)$ when $ d=6$). 

We observe that, if we do not make further hypothesis on the positive solution $u$, the non-negativity assumption on the Ricci tensor is necessary to obtain the rigidity result of the underlying manifold. Indeed, there are hyperbolic-type model manifolds supporting bounded, positive and radially symmetric solutions to \eqref{eq_crit} (see \cite{BFG}, \cite{MS}). Nevertheless, the rigidity result still holds true if one assumes that $M$ is a Cartan–Hadamard model manifold and the solution  to \eqref{eq_crit} is  radially symmetric and of finite energy or that it minimizes the Sobolev quotient \cite{MS}.

\medskip

In the case of Riemannian surfaces we have the following new and optimal result 

\begin{theorem} \label{thm_Liouville}
Let $(M,g)$ be a complete, connected, non-compact, boundaryless Riemannian manifold 
of dimension $2$ and with nonnegative Ricci curvature.  
Let $G=G(t)$ be any positive nondecreasing function satisfying $ \int_c^{\infty} \frac{dt}{tG(t)} = + \infty,$ for some $c>1$ and denotes by $r=r(x)$ the geodesic distance from a fixed reference point of $M$. 

If $u \in C^3(M)$ is a solution to 
\begin{equation} \label{eq_crit2}
-\Delta u = e^u  \qquad \textmd{ in } \quad M, 
\end{equation} 
such that
\begin{equation}\label{eq_crit2-hyp}
u(x) \geq - 4 \log \left[ r(x)G^{\frac{1}{2}}(r(x))\right], \qquad \textmd{for} \quad  r(x) >c,
\end{equation}
then $(M,g)$  is isometric to $\R^2$ with the Euclidean metric and $u$ is given by
%Then $(M,g)=(\R^2,g_E)$ and $u$ is given by
\begin{equation}\label{u_radial_liouv}
u(x) = \log \frac{1}{(a+b|x-x_0|^2)^2} \,,
\end{equation}
for some $a,b>0$ with $1=8ab$ and $ x_0 \in \R^2$. 

Furthermore, the coefficient $4$ in the asymptotic lower bound \eqref{eq_crit2-hyp} is optimal. Indeed, for any $\beta >4$ there is a non flat, complete, connected, non-compact, boundaryless Riemannian manifold $(M,g)$ of dimension $2$, with nonnegative Ricci curvature and supporting a  (finite mass) solution to the Liouville equation \eqref{eq_crit2} satisfying $u(x) \geq - \beta \log \left[ r(x) \right],$ for $ r(x) >> 1$. 

\end{theorem}

\medskip

Note that the above result is new even in the Euclidean case. Indeed, differently from all the previous works, we do not assume the finite mass condition \eqref{finite-mass-Liouv_Rd}. 

We also notice that the lower bound \eqref{eq_crit2-hyp} is quite natural since, in the Euclidean case, it is satisfied by all the solutions having finite mass (just choose $G(t)$ a positive constant larger than $b$ in \eqref{eq_crit2-hyp} and $c$ large enough, depending on $b$). 

\smallskip

We also observe that Theorem 1.6 in \cite{CatinoMonticelli} is a very particular case of our Theorem \ref{thm_Liouville}. Indeed, in \cite{CatinoMonticelli} the authors assume both the finite mass condition 
\begin{equation} \label{finite-mass-Liouv_M}
\int_{M} e^u\, <+\infty \,
\end{equation} 
and the lower bound \eqref {eq_crit2-hyp} with the special function $G(t)= \log^{\gamma}(t)$, 
$ \gamma \in  [0,1)$.\footnote{Observe that, for any $k \geq1$, $G_k(t) = \prod_{j=1}^k L_j(t)$, where $L_1(t)=\log(t)$ and $L_j(t) = \log(L_{j-1}(t))$ if $ j \geq 2$, are (more general) functions satisfying the assumption of Theorem \ref{thm_Liouville}. }

\smallskip

Finally we mention the very recent work \cite{CaiLai} where the conclusion of  Theorem  \ref{thm_Liouville} was obtained by assuming both the finite mass condition \eqref{finite-mass-Liouv_M} and the asymptotic lower bound : $u(x) \geq - 4 \log \left[ r(x)\right] + o(\log \left[ r(x)\right] ),$ as $r(x) \longrightarrow +\infty.$

\smallskip

As already noted, the lower bound \eqref {eq_crit2-hyp} of Theorem \ref{thm_Liouville} is optimal.
Nevertheless, the flexibility of our approach allows us to weaken the lower bound \eqref {eq_crit2-hyp} and still achieve a (new) classification result. This is the content of the next result.

\begin{theorem} \label{thm_Liouville2}
Let $(M,g)$ be a complete, connected, non-compact, boundaryless Riemannian manifold 
of dimension $2$ and with nonnegative Ricci curvature.  
Let $G=G(t)$ be any positive nondecreasing function satisfying $ \int_c^{\infty} \frac{dt}{tG(t)} = + \infty,$ for some $c>1$ and denotes by $r=r(x)$ the geodesic distance from a fixed reference point of $M$. 

If $u \in C^3(M)$ is a solution to 
\begin{equation} \label{eq_crit3}
-\Delta u = e^u  \qquad \textmd{ in } \quad M, 
\end{equation} 
such that
\begin{equation}\label{eq_crit3-hyp}
u(x) \geq - 6 \log \left[ r(x)G^{\frac{1}{3}}(r(x))\right], \qquad \textmd{for} \quad  r(x) >c,
\end{equation}
then $(M,g)$  is conformal to $\R^2$ with the Euclidean metric.

Furthermore, the result is optimal in the sense that there are complete manifolds $(M,g)$ of dimension $2$, with nonnegative Ricci curvature, conformal to $\R^2$ with the Euclidean metric, and supporting a (finite mass) solution to the Liouville equation \eqref{eq_crit2} satisfying the lower bound \eqref{eq_crit3-hyp}. 
\end{theorem}

Once again, the previous result is obtained without making the assumption of finite mass \eqref{finite-mass-Liouv_M}.

\bigskip

In the next two results we study non-negative  solutions to subcritical equations. In order to fix the notations, we set $p_S$ the critical exponent, i.e., 
%We will denote by $p_S$ the critical exponent appearing in \eqref{eq_crit}:
\begin{equation} \label{p_S}
p_S= 
\begin{cases}
\dfrac{d+2}{d-2} & \textmd{ if } d \geq 3 \\
+\infty & \textmd{ if } d=2 \,.
\end{cases}
\end{equation}

Then,  we have

\begin{theorem} \label{thm_subcritical}
Let $(M,g)$ be a complete, connected, non-compact, boundaryless Riemannian manifold 
of dimension $d\geq 2$ and with nonnegative Ricci curvature.  Let $u \in C^3(M)$ be a nonnegative solution to 
\begin{equation} \label{eq_subcrit}
-\Delta u = u^{p}  \qquad \textmd{ in } \quad M,
\end{equation}
with $1<p<p_S$.  
Then $u=0$.
\end{theorem} 

Theorem \ref{thm_subcritical} was already known in literature, and it goes back to \cite{GS} (see also \cite{CatinoMonticelli}). Hence, the novelty of the theorem is not in the result itself but the way we obtain it. Indeed, our approach appears to be very simple and for this reason we think that it deserves to be emphasized.  

\medskip

In the subcritical case, we can also weaken the condition on the Ricci curvature whenever the manifold $M$ satisfies a Sobolev inequality of the form
\begin{equation} \label{Sobolev_nu}
\|\nabla \psi \|_{L^2(M)} \geq S \|\psi\|_{L^{\nu}(M)} \qquad \forall \,  \psi \in C_c^1(M),
\end{equation}
for some $\nu >2$ and where $d\geq 3$.

More precisely, we define $\Ric_{-} : M \mapsto \left[0, \infty \right)$ by $\Ric_{-}(x) = 0$ if $\Ric(x) \geq 0$ and, if 
$\Ric(x)$ has a negative eigenvalue, then $- \Ric_{-}(x)$ is the lowest eigenvalue of $\Ric(x).$ 

\smallskip

In this context, our new result is the following   

\begin{theorem} \label{thm_subcritical_Ricci_neg}
Let $(M,g)$ be a complete, connected, non-compact, boundaryless Riemannian manifold 
of dimension $d\geq 3$ and let $u \in C^3(M)$ be a nonnegative solution to 
\begin{equation} \label{eq_subcrit}
-\Delta u = u^{p}  \qquad \textmd{ in } \quad M,
\end{equation}
with $1<p<p_S$.

\noindent Let $M$ enjoy a Sobolev inequality \eqref{Sobolev_nu} for some $\nu >2$ and assume that 
\begin{equation} \label{Ricci_bound_neg}
\|\Ric_{-}\|_{L^{\frac{\nu}{\nu-2}}(M)}  \leq \frac{1}{48}(p_S-p)(p-1)(d-2) S^2  \,.
\end{equation}
Assume also that, for a fixed point $x_0 \in M$,  the volume of the geodesic ball $B_R= B_R(x_0)$ satisfies 
\begin{equation} \label{Vol_bound_in_p}
V_g(B_R) = O( R^{2 + \frac{8}{p-1}}) \qquad  \textmd{as}  \quad R\to \infty. 
\end{equation}
Then $u=0$. 
\end{theorem} 

We observe that, for any $d\geq 3$ and any $1<p<p_S$, there is a complete, connected, non-compact, boundaryless Riemannian manifold with non trivial $\Ric_{-}$ and satisfying all the assumptions of Theorem \ref{thm_subcritical_Ricci_neg}.  See item (vi) of Remark  \ref{remark_intro}. 

\medskip

There are some further noteworthy observations on Theorem \ref{thm_subcritical_Ricci_neg}, which we collect in the following remark.

\begin{remark} \label{remark_intro} 
$\ $

\begin{itemize}
\item[(i)] It  is interesting to note that, in the two-dimensional case,  the conclusion of Theorem \ref{thm_subcritical_Ricci_neg} holds true under the sole assumption that the manifold $M$
has polynomial volume growth (of any order). See Theorem \ref{thm_serrin}.

\item[(ii)] The exponent appearing in  the volume growth condition  \eqref{Vol_bound_in_p} is always larger than $d$, the dimension of the manifold $M.$ Indeed, $2 + \frac{8}{p-1} > 2(d-1) >d,$ since $d \geq 3$ and  $1<p<p_S$.
\item[(iii)] As we will explain in Section \ref{section2}, it will be useful to introduce the parameter 
\begin{equation} \label{n_def_intro}
n= 2\, \dfrac{p+1}{p-1} \,.
\end{equation}
Notice that  $n \geq d$, and $n=d$ if and only if $p=p_S$. Hence, in terms of $n$, conditions \eqref{Ricci_bound_neg} and \eqref{Vol_bound_in_p} can be written as 
\begin{equation} \label{K-cond}
\|\Ric_{-}\|_{L^{\frac{\nu}{\nu-2}}(M)}  \leq   \frac{n-d}{3(n-2)^2} S^2 
\end{equation}
and
\begin{equation} \label{Vol_bound}
V_g(B_R) = O( R^{2(n-1)}) \quad \textmd{as } R \to \infty\,,
\end{equation}
respectively.

\begin{equation}% \label{Vol_bound}
\end{equation}

\item[(iv)] When $M$ enjoys an Euclidean type Sobolev inequality, i.e., when $M$ satisfies  
\eqref{Sobolev_nu} with $\nu=\dfrac{2d}{d-2} $,  then \eqref{Ricci_bound_neg} boils down to 
\begin{equation} \label{Ricci_bound_neg-euclideo}
\|\Ric_{-}\|_{L^{\frac{d}{2}}(M)}  \leq  \frac{1}{48}(p_S-p)(p-1)(d-2) S^2 \,.
\end{equation}

\item[(v)] The volume growth in \eqref{Vol_bound_in_p} is automatically satisfied  in many situations. 
For instance, if one assumes that $\Ric_- \in L^{\frac{d}{2}}(M) \cap  L^{\frac{\mu}{2}}(M)$, for some $\mu>d$, then there exists a positive constant  $R_0=R_0(d,\Ric_-)$ such that
$$
V_g(B_R)  \leq C(d) R^d \left(\log\left( \frac{2R}{R_0} \right) \right)^{\frac{d}{2} - 1} \quad \textmd{if } R \geq R_0 \,,
$$
(see for instance the survey \cite{Carron} and also \cite{Aubry, Gallot, LiYau, PetersenWei, PigolaRigoliSetti} for other related results). Hence, thanks to Remark (ii), \eqref{Vol_bound_in_p} is satisfied.

\item[(vi)] Hereafter we provide an example of a complete, connected, non-compact, boundaryless Riemannian manifold with non trivial $\Ric_{-}$ and satisfying all the assumptions of Theorem \ref{thm_subcritical_Ricci_neg}.  The example can be seen as a ``small perturbation"  of the Euclidean space $\R^d,$ $d \geq 3$.  
We are indebted to  Matteo Muratori, who suggested this example.

For any $d\geq 3$ and any $1<p<p_S$ consider the model manifold $(M,g)$, namely the manifold with a pole $o$ and whose metric is given, in polar coordinates around $o$, by 

\begin{equation}% \label{Vol_bound}
ds^2 = dr^2 + \psi^2(r)d\theta^2 \qquad r>0, \quad \theta \in \mathbb{S}^{d-1},
\end{equation}
where $d\theta^2$ is the standard Euclidean metric on $\mathbb{S}^{d-1}$ and 
\begin{equation}% \label{Vol_bound}
\psi(r) = r + \varepsilon \eta(r), \qquad r>0,
\end{equation}
with $ \varepsilon \in (0,1)$ and $ \eta$ a fixed smooth cut-off function on $(0, +\infty)$ such that : $ \eta(r) = 0 $ for $r \in (0,\frac{1}{2}) \bigcup (3,+\infty)$,  $ \eta(r) = 1 $ for $r \in (1,2) $ and  $ 0 \leq \eta(r) \leq 1$ for $r >0.$  

By definition of $\psi$, for any $ \varepsilon \in (0,1)$  we have 
\begin{equation} \label{doppio-bound-psi}
r \leq \psi(r) \leq 4r \qquad \textmd{ for } r>0,
\end{equation}
hence there is a constant  $C_1(d)>0$, depending only on $d$, such that $V_g(B_R(o))  \leq C_1(d) R^d$ for any $ R>0$. Thus \eqref{Vol_bound_in_p} is satisfied. 

Also, by construction there exists a constant $C_2(d,\eta)>0$, depending only on $d$ and $\eta$, such that for any 
$ \varepsilon \in (0,1)$
\begin{equation*} 
\begin{cases}
\Ric_{-} \not \equiv 0 & \textmd{ on } M, \\
\Ric_{-} \leq \varepsilon C_2(d, \eta) & \textmd{ on } M, \\
\Ric_{-} = 0 &  \textmd{ on } M \setminus B_4(o), 
\end{cases}
\end{equation*}
and so 
\begin{equation}\label{Ricci-modello}
0 < \|\Ric_{-}\|_{L^{\frac{d}{2}}(M)}  \leq  \varepsilon C_2(d,\eta) \left[V_g(B_4(o))\right]^{\frac{2}{d}} \leq  \varepsilon C_2(d,\eta)
C_1(d)^{\frac{2}{d}} 4^2.
\end{equation}

Since \eqref{doppio-bound-psi} is in force, we see that,  for any $\varphi  \in C^1_c(M)$ and for any 
$ \varepsilon \in (0,1)$
\begin{equation*} 
\begin{cases}
\| \varphi \|_{L^{\frac{2d}{d-2}}(\R^d)}  \leq  \| \varphi \|_{L^{\frac{2d}{d-2}}(M)} \leq C_3(d) \| \varphi \|_{L^{\frac{2d}{d-2}}(\R^d)}, \\
\| \nabla \varphi \|_{L^2 (\R^d)}  \leq  \| \nabla \varphi \|_{L^{2}(M)} \leq C_3(d) \|  \nabla \varphi \|_{L^{2}(\R^d)},
\end{cases}
\end{equation*}
for some constant $C_3(d)>1$, depending only on the dimension $d$. 
Therefore, $M$ enjoys the following Euclidean type Sobolev inequality,
\begin{equation*} 
\|\nabla \varphi \|_{L^2(M)} \geq S(\R^d) \left[ C_3(d) \right]^{-1} \|\varphi\|_{L^{\frac{2d}{d-2}}(M)} 
\qquad \forall \,  \varphi \in C_c^1(M), \qquad \forall \, \varepsilon \in (0,1),
\end{equation*}
where $S(\R^d)$ denotes the best constant in the Sobolev inequality in the  Euclidean space $\R^d$. Since the positive constant $S(\R^d) \left[ C_3(d) \right]^{-1}$ is independent of $\varepsilon \in (0,1)$, we see that 
also \eqref{Ricci_bound_neg-euclideo} is satisfied by taking $\varepsilon$ small enough in \eqref{Ricci-modello}.

\end{itemize}
\end{remark}

\medskip

\noindent {\bf Strategy of the proof.} The proofs of our main results 
%Theorems \ref{thm_critical}, \ref{thm_subcritical}, \ref{thm_subcritical_Ricci_neg} and \ref{thm_Liouville} 
exploit a common strategy which is related to the definition of a suitable $P$-function. Indeed, after letting
$$
u=v^{-\frac{2}{p-1}}
$$
for equations \eqref{eq_crit}, where $p=p_S$, and \eqref{eq_subcrit}, and 
$$
u=-2 \log v
$$
for \eqref{eq_crit2}, we find that $v$ satisfies
$$
\Delta v = P \,,
$$
where
$$
P = \frac{1}{v} \left(\frac{n}{2} |\nabla v|^2 + \frac{2}{n-2} \right) 
$$
with
$$
\frac{n}{2} = \frac{p+1}{p-1} \,,
$$
(notice that $n \geq d$ and that $n=d$ when $p=p_S$) and 
$$
P = \frac{1}{v} \left( |\nabla v|^2 + \frac12 \right) 
$$
for $d=n=2$. The main goal in Theorems \ref{thm_critical} and \ref{thm_Liouville} is to prove that $P$ is constant, which implies the classification results, while Theorems \ref{thm_subcritical} and \ref{thm_subcritical_Ricci_neg} follow by proving a contradiction argument on $P$. In both cases, this is done by noticing that $P$ satisfies
$$
\diver(v^{2-n} \nabla P ) =  nv^{1-n}k[v] \,,
$$
where 
$$
k[v] = \Big{|} \nabla^2 v - \frac{\Delta v}{d} \Big{|}^2 + \frac{n-d}{nd} (\Delta v)^2 + \Ric(\nabla v, \nabla v) \,.
$$
By exploiting the definitions of $P$ and of $k[v]$, we obtain that 
\begin{equation*} 
\divergence \left( P^{t-1}  v^{2-n}\nabla P \right)  \geq \left(t-\frac12\right) P^{t-2} v^{2-n}|\nabla P|^2 + nP^{t-1}v^{1-n}  \left\{ \Ric(\nabla v,\nabla v) + \frac{n-d}{n^2} (\Delta v)^2 \right\}  
\end{equation*}
for any $t \in \R$. Since $\Delta v = P$, by letting $\psi_R$ a suitable cut-off function, one easily obtains that
\begin{multline*}
 \left(t-\frac12\right) \int_{B_R} P^{t-2} v^{2-n} |\nabla P|^2  
 %{\color{red}{- }} n \int_{B_R} P^{t-1}v^{1-n}   {\color{red}{\Ric_-|\nabla v|^2}}
%\Ric_-(\nabla v,\nabla v) 
+ \frac{n-d}{n} \int_{B_R}  P^{t+1}v^{1-n}   \\
\leq n \int_{B_R} P^{t-1}v^{1-n}   \Ric_-|\nabla v|^2 - \int_{B_{2R} \setminus B_{R}}  P^{t-1}  v^{2-n}\nabla P \cdot \nabla \psi_R \,.
\end{multline*}

Then the conclusion follows once we can provide suitable bounds for the RHS in the inequality above for $t> 1/2$. This will be done by carefully choosing the parameter $t$, as well as other parameters appearing in the proofs, and by letting $R \to \infty$. 
%Crucial information will come from the volume growth of geodesic balls in $M$ and by using Bishop-Gromov Theorem. 
In the case of Theorems \ref{thm_subcritical} and \ref{thm_subcritical_Ricci_neg} these estimates will force a contradiction (unless $u\equiv 0$ in $M$), while in the case of Theorems \ref{thm_critical} and \ref{thm_Liouville} the conclusion will follow by showing that $P$ is constant, implying that $k[v]=0$, and a classical splitting theorem (see \cite[Theorem 5.7.4]{Petersen} and also \cite{FMM} and \cite{CatinoMonticelli}) yields the classification result.

\medskip

\noindent {\bf  Organization of the paper.} In Section \ref{section2} we introduce the $P$-function and show that it is a subsolution of a suitable equation. In Section \ref{section3} we prove the results for subcritical equations Theorems \ref{thm_subcritical} and \ref{thm_subcritical_Ricci_neg}. We prefer to start from these theorems since, in this case, the results follow almost straightforwardly from the discussion above and by arguing by contradiction.  Then we move to the classification results in Theorems \ref{thm_critical} and \ref{thm_Liouville}. In Section \ref{section_integral_estimates} we prove some preliminary integral estimates which are particularly useful for the proofs of Theorems \ref{thm_critical} and \ref{thm_Liouville}. Then in Sections \ref{section5} and \ref{Section_Liouville} we will give the proof of Theorems \ref{thm_critical} and \ref{thm_Liouville}, respectively, by showing that $P$ is constant.

\section{The p-function} \label{section2}
In this section we introduce an auxiliary function, the so-called $P$-function and prove some preliminary results.
Let $u$ be a solution to 
\begin{equation} \label{eq_semilinear1}
- \Delta u = u^p \qquad \textmd{ in } \quad M,
\end{equation}
where $1<p \leq p_S$ when $d\geq 3$, and $1<p<+\infty$ for $d=2$. 
Let 
\begin{equation} \label{v_def}
u=v^{-\frac{2}{p-1}} \,.
\end{equation}
Notice that if $p=p_S$ then \eqref{v_def} reads as $u=v^{-\frac{d-2}{2}}$. Since
$$
\nabla u = -\frac{2}{p-1} v^{-\frac{2}{p-1}-1} \nabla v \quad \textmd{ and } \quad \Delta u = 2 \frac{p+1}{(p-1)^2} v^{-\frac{2}{p-1} -2 }|\nabla v|^2 - \frac{2}{p-1} v^{-\frac{2}{p-1}-1} \Delta v
$$
it is readily seen that $v$ satisfies
$$
\Delta v = \frac{p+1}{p-1} \frac{|\nabla v|^2}{v} + \frac{p-1}{2} \frac{1}{v}\,,
$$
which we write as
\begin{equation} \label{eq_v}
\Delta v = P \,,
\end{equation}
where
\begin{equation}\label{P_def1}
P = \frac{1}{v} \left(\frac{n}{2} |\nabla v|^2 + \frac{2}{n-2} \right) 
\end{equation}
with
$$
\frac{n}{2} = \frac{p+1}{p-1} \,,
$$
i.e.
$$
p = \frac{n+2}{n-2} \,.
$$
We notice that $n \geq d$ and $n=d$ for $p=p_S$ and $d\geq3$. If $d=2$ then $n>d$.  

In the case $d=2$, we also consider the Liouville equation 
\begin{equation} \label{eq_semilinear2}
- \Delta u = e^u \qquad \textmd{ in } \quad M,
\end{equation}
and, by letting 
\begin{equation} \label{v_def2}
u=-2 \log v
\end{equation}
we see that $v$ satisfies \eqref{eq_v} where now 
\begin{equation}\label{P_def2}
P = \frac{1}{v} \left(|\nabla v|^2 + \frac12 \right) \,,
\end{equation}
and in this case we set $n=d=2$. 

This is our starting point. Indeed, by letting
\begin{equation} \label{n_def}
n  = 
\begin{cases}
2 \dfrac{p+1}{p-1} & \textmd{ if $u$ is a solution to \eqref{eq_semilinear1}} \,,\\
2 &  \textmd{ if $u$ is a solution to \eqref{eq_semilinear2}} \,,
\end{cases}
\end{equation}
and noticing that $n \geq d$, we have that both \eqref{eq_semilinear1} and \eqref{eq_semilinear2} can be reduced to the study of positive solutions to 
$$
\Delta v = P \,,
$$
where
\begin{equation}\label{P_def12}
P = \dfrac{1}{v} \Big(\dfrac{n}{2} |\nabla v|^2 + c_n \Big)
\end{equation}
with
\begin{equation}\label{cn_def}
c_n= 
\begin{cases}
\dfrac{2}{n-2} & \textmd{ if } n>2 \\
& \\ 
\dfrac12  & \textmd{ if } n=2 \,.
\end{cases}
\end{equation}
The right-hand side of \eqref{eq_v}, which is defined in \eqref{P_def12}, is what we call $P$-function. Our main goal in the proof of Theorems \ref{thm_critical} and \ref{thm_Liouville} will be to show that $P$ is constant, while we will prove nonexistence results by contradiction in the case of Theorems \ref{thm_subcritical} and \ref{thm_subcritical_Ricci_neg}. The crucial property of $P$ is that $P$ is a positive subsolution to a suitable elliptic equation, as it is proved in the following lemmas.

Before stating these lemmas, it will be useful to introduce the following quantity 
\begin{equation} \label{k_def1} 
k[v]=| \nabla^2 v|^2 -\frac{1}{n} (\Delta v)^2+ \Ric(\nabla v,\nabla v) \,.
\end{equation}
It is straightforward to verify that
\begin{equation} \label{k_def2} 
k[v]=\Big{|} \nabla^2 v - \frac{\Delta v}{d} g \Big{|}^2 +\frac{n-d}{nd} (\Delta v)^2 + \Ric(\nabla v,\nabla v)  \,.
\end{equation}
Notice that, since $n \geq d$, if $\Ric \geq 0$ then $k[v] \geq 0$.

In the first lemma we prove a differential identity which is crucial for our approach. 

\begin{lemma} \label{lemma_appa}
Let $(M,g)$ be a Riemannian manifold of dimension $d\geq 2$. Let $P$ and $n$ be given by \eqref{P_def12} and \eqref{n_def}, respectively, and let $v \in C^3(M)$, with $v>0$. Then we have
\begin{equation*}
\divergence \bigg( v^{1-n}\bigg( \frac{n}{2}\nabla |\nabla v|^2 - n \Delta v \nabla v + (n-1)P \nabla v\bigg)\bigg) \\
 =n v^{1-n}k[v] + v \bigg(\Delta v- P \bigg) \Delta(v^{1-n}),
\end{equation*}
where $k[v]$ is given by \eqref{k_def1}.
\end{lemma}

\begin{proof}
We first notice that
\begin{equation*}
\frac{n}{2} \divergence \left(  v^{1-n}  \nabla |\nabla v|^2  \right)  = \frac{n}{2} v^{1-n} \Delta \left( |\nabla v|^2   \right) + (1-n) \frac{n}{2} v^{-n} \nabla v \cdot \nabla |\nabla v|^2  \,,
\end{equation*}
and Bochner identity yields
\begin{equation} \label{appa1.1}
\frac{n}{2}\divergence \left( v^{1-n}  \nabla |\nabla v|^2  \right)  = n v^{1-n}  \left( |\nabla^2 v|^2 + \nabla v \cdot \nabla \Delta v  + \Ric( \nabla v ,\nabla v)   \right)  + (1-n) \frac{n}{2} v^{-n} \nabla v \cdot \nabla |\nabla v|^2 \,.
\end{equation}
Since
\begin{equation*}
-n  \divergence \left(  v^{1-n}  \Delta v\nabla v  \right)  = -n \left[ v^{1-n} (\Delta v)^2 + v^{1-n} \nabla v \cdot \nabla \Delta v + (1-n)v^{-n} \Delta v |\nabla v|^2 \right] \,,  \end{equation*}
from \eqref{appa1.1} we obtain 
\begin{multline*}
\divergence \left( \frac{n}{2} v^{1-n}  \nabla |\nabla v|^2  -nv^{1-n}  \Delta v\nabla v  \right) 
= n v^{1-n}  \left[ |\nabla^2 v|^2 - (\Delta v)^2   + Ric( \nabla v ,\nabla v)   \right] \\
+ (1-n) v^{-n} \left(  \frac{n}{2} \nabla|\nabla v|^2 \cdot \nabla v  -n \Delta v |\nabla v|^2   \right) \,,
\end{multline*}
i.e. 
\begin{multline} \label{appa1.2}
 \divergence \left\{  \left( \frac{n}{2} v^{1-n}  \nabla |\nabla v|^2  -nv^{1-n}  \Delta v\nabla v  \right) \right\}
= n v^{1-n} k[v] + (1-n)  v^{1-n}   (\Delta v)^2  \\
+ (1-n) v^{-n} \left(  \frac{n}{2} \nabla|\nabla v|^2 \cdot \nabla v  -n \Delta v |\nabla v|^2   \right) \,.
\end{multline}
%and, since $ \frac{n}{2}\nabla  |\nabla v|^2 = \nabla(vP) = P \nabla v + v \nabla P$ and by the definition of $k[v]$, we obtain
%\begin{multline*}
%\Cl \left( \frac{n}{2} v^{1-n}  \nabla |\nabla v|^2  -nv^{1-n}  \Cl v\nabla v  \right) 
%= n v^{1-n} k[v] + (1-n)  v^{1-n}   (\Cl v)^2 + \\
%+ (1-n) v^{-n} \left( P |\nabla v |^2   -n \Cl v |\nabla v|^2   \right)  + (1-n) v^{1-n} \nabla P \cdot \nabla v \,.
%\end{multline*}
%
%
Now we consider
\begin{equation*}
\begin{split}
(n-1) \divergence & (v^{1-n} P \nabla v) \\
&  = (n-1) v^{1-n} \left( P \Delta v  + \nabla P \cdot \nabla v \right) - (n-1)^2 v^{-n}  P |\nabla v|^2  \\
& =  (n-1) v^{1-n} \left( P \Delta v  - \frac{P}{v} |\nabla v|^2 + \frac{n}{2v} \nabla |\nabla v|^2 \cdot \nabla v \right) - (n-1)^2 v^{-n}  P |\nabla v|^2\,,
\end{split}
\end{equation*}
where in the last equality we used  
$$
\nabla P = - \frac{P}{v} \nabla v + \frac{n}{2v} \nabla |\nabla v|^2\,,
$$
and hence
\begin{equation} \label{appa1.3}
(n-1) \divergence ( v^{1-n} P \nabla v) \\ =  (n-1) v^{1-n}  P \Delta v   +  \frac{n}{2}  (n-1) v^{-n} \nabla |\nabla v|^2 \cdot \nabla v   - n(n-1) v^{-n}  P |\nabla v|^2 \,.
\end{equation}
Thus from \eqref{appa1.2} and \eqref{appa1.3} we find
\begin{multline*}
\divergence \left(   v^{1-n} \left[ \frac{n}{2} \nabla |\nabla v|^2  -n  \Delta v\nabla v + (n-1)  P \nabla v \right] \right) 
= n v^{1-n} k[v] + (1-n)  v^{1-n}   (\Delta v)^2  \\
+ n (n-1) v^{-n}  \Delta v |\nabla v|^2 + (n-1) v^{1-n} P   \Delta v  - n (n-1) v^{-n} P  |\nabla v|^2   \,,
 \end{multline*}
i.e.
\begin{multline*}
\divergence \left( v^{1-n} \left[ \frac{n}{2} \nabla |\nabla v|^2  -n   \Delta v\nabla v + (n-1)   P \nabla v \right] \right) \\
= n v^{1-n} k[v]  + n(n-1) v^{-n}   |\nabla v|^2 (\Delta v - P) + (n-1) v^{1-n}  \Delta v (P- \Delta v)
\end{multline*}
and thus
\begin{equation*}
\divergence \left(  v^{1-n} \left[ \frac{n}{2} \nabla |\nabla v|^2  -n   \Delta v\nabla v + (n-1)   P \nabla v \right] \right) \\
= n v^{1-n} k[v]  + (n-1) v^{-n} (n  |\nabla v|^2 - v \Delta v)(\Delta v - P)   \,.
\end{equation*}
Since
$$
v \Delta (v^{1-n}) = (n-1) v^{-n} \left[ n  |\nabla v|^2 - v \Delta v \right]
$$
we conclude.
\end{proof}

From Lemma \ref{lemma_appa} we readily obtain the following crucial differential identity for the $P$-function.

\begin{lemma} \label{lemma_P}
Let $(M,g)$ be a Riemannian manifold of dimension $d\geq 2$. Let $n$ and $P$ be given by \eqref{n_def} and \eqref{P_def12}, respectively. If $v$ is a positive solution to \eqref{eq_v} then $P$ satisfies
\begin{equation}\label{eqfond}
\divergence \left(v^{2-n}\nabla P\right)=nv^{1-n}k[v] \quad \textmd{ in } M \,.
\end{equation}
\end{lemma}

\begin{proof}
From Lemma \ref{lemma_appa} and since $\Delta v = P$ we have
\begin{equation} \label{2}
\divergence \left( v^{1-n}\left( \frac{n}{2}\nabla |\nabla v|^2 - \Delta v \nabla v \right)\right)   =n v^{1-n}k[v] \,.
\end{equation}
Since $\Delta v =P$ and from \eqref{P_def12}, we have that
\begin{equation} \label{nablaP}
\nabla P =   \frac{n}{2v}\nabla |\nabla v|^2  -\frac{1}{v} P \nabla v = \frac{1}{v} \left( \frac{n}{2}\nabla|\nabla v|^2 -  \Delta v \nabla v \right)
\end{equation}
which, together with \eqref{2}, implies \eqref{eqfond}.
\end{proof}

Lemma \ref{lemma_P} is also preliminary to the following lemma, which provides another fundamental differential inequality. 
%We mention that, when $M$ is the Euclidean space $\R^d$ and $n=d$, an inequality similar to the one of the following lemma can be found in the proof of the main results of \cite{OU} and \cite{Vetois}. Indeed, this result is the starting point in \cite{OU} and \cite{Vetois} for classifying solutions to the critical $p$-Laplace equation in $\R^n$ without energy assumptions.

\begin{lemma} \label{lemma_Pt}
Let $(M,g)$ be a Riemannian manifold of dimension $d\geq 2$. Let $n$ and $P$ be given by \eqref{n_def} and \eqref{P_def12}, respectively. If $v$ is a positive solution to \eqref{eq_v} then
\begin{equation} \label{eqPt}
\divergence \left( P^{t-1}  v^{2-n}\nabla P \right)  \geq \left(t-\frac12\right) P^{t-2} v^{2-n}|\nabla P|^2 + nP^{t-1}v^{1-n}  \left\{ \Ric(\nabla v,\nabla v) + \frac{n-d}{n^2} P^2 \right\}  
\end{equation}
in $M$, for any $t \in \R$.
\end{lemma}

\begin{proof}
We first notice that from \eqref{nablaP}, we have that 
$$
\nabla P=\frac{n}{v}\left( \nabla^2  v-\frac{\Delta v}{n} g \right) \nabla v \,,
$$
and Cauchy-Schwarz inequality yields
\begin{equation*}
|\nabla P|^2\le\frac{n^2}{v^2}\Big{|}\nabla^2 v-\frac{\Delta v}{n}{g}\Big{|}^2|\nabla v|^2 \,.
\end{equation*}
Since $P>0$, for any $t\in\R$ we can write
\begin{equation}\label{comp1}
P^{t-2}v^{2-n}|\nabla P|^2\le P^{t-2}v^{2-n}\frac{n^2}{v^2}|\nabla^2 v-\frac{\Delta v}{n}{g}|^2|\nabla v|^2= P^{t-1}v^{1-n}n^2|\nabla^2 v-\frac{\Delta v}{n} g|^2\frac{|\nabla v|^2}{vP}.
\end{equation}
From the definition of $P$ \eqref{P_def12} we have
\begin{equation}\label{comp2}
\frac{|\nabla v|^2}{vP}\le\frac{2}{n},
\end{equation}
and then \eqref{comp1} and \eqref{comp2} yield
\beq\label{comp3}
P^{t-2}v^{2-n}|\nabla P|^2\le P^{t-1}v^{1-n}2n\Big{|}\nabla^2 v-\frac{\Delta v}{n} g \Big{|}^2 \,.
\eeq
Now we notice that
\begin{equation*}
\Big{|}\nabla^2 v-\frac{\Delta v}{n} g\Big{|}^2 =|\nabla^2 v|^2 + \frac{d}{n^2} (\Delta v)^2 - \frac{2}{n}  (\Delta v)^2 = |\nabla^2 v|^2 - \frac{1}{n}  (\Delta v)^2 
 - \frac{n-d}{n^2} (\Delta v)^2 
 \end{equation*}
and hence from \eqref{k_def1} we find
\begin{equation*}
\Big{|}\nabla^2 v-\frac{\Delta v}{n} g\Big{|}^2 =k[v]-  \Ric(\nabla v,\nabla v) - \frac{n-d}{n^2} (\Delta v)^2 
 \end{equation*}
and from \eqref{comp3} we have
\begin{equation*} 
P^{t-2}v^{2-n}|\nabla P|^2\le P^{t-1}v^{1-n}2n \left\{ k[v] -  \Ric(\nabla v,\nabla v) - \frac{n-d}{n^2} (\Delta v)^2 
 \right\} \,,
\end{equation*}
i.e.
\begin{equation} \label{starb}
P^{t-2}v^{2-n}|\nabla P|^2 + 2nP^{t-1}v^{1-n}  \left\{ \Ric(\nabla v,\nabla v) + \frac{n-d}{n^2} (\Delta v)^2 
 \right\} \le 2n P^{t-1}v^{1-n} k[v] \,.
\end{equation}
Since 
\begin{equation*}
\divergence \left(P^{t-1}  v^{2-n}\nabla P \right)  =P^{t-1} \diver \left( v^{2-n}\nabla P \right) + (t-1) P^{t-2} v^{2-n}|\nabla P|^2 \,,
\end{equation*}
from \eqref{eq_v}, \eqref{eqfond} and \eqref{starb} we obtain \eqref{eqPt}.
\end{proof}

\section{Classification results for subcritical equations} \label{section3}
In this section we prove Theorems \ref{thm_subcritical} and  \ref{thm_subcritical_Ricci_neg}. We notice that, in this case, we only use Lemma \ref{lemma_Pt} and an information on the volume growth of geodesic balls.

We first consider the case of manifolds with nonnegative Ricci curvature. We recall that, since $\Ric \geq 0$, Bishop-Gromov relative volume comparison theorem \cite[Lemma 7.1.4]{Petersen} implies 
\begin{equation} \label{BG}
V_g(B_r) \leq \omega_d r^d, \qquad \forall \, r>0,
\end{equation}
where $\omega_d$ is the Euclidean volume of the ball.

To proceed further we recall that, having fixed a point $x_0 \in M$,  for every $R>0$ there exists a (standard) cutoff function $\varphi=\varphi_R$, i.e., a function $\varphi=\varphi_R \in C^{0,1}_c(M)$ such that : 
\begin{equation}\label{vphi_def}
\varphi =1 \textmd{ in } B_{R}, \quad \varphi = 0 \textmd{ in } M \setminus  B_{2R} \qquad \textmd{ and  }
\qquad 0 \leq \varphi \leq 1, \quad \vert \nabla \varphi \vert \leq \frac{1}{R}\quad  \textmd{ on } M,
\end{equation}
where $B_r = B_r(x_0)$ denotes the geodesic ball of center $x_0$ and radius $r>0$.

\begin{proof}[Proof of Theorem \ref{thm_subcritical}]
When $d=2$, the volume growth \eqref{BG} implies that $M$ is a parabolic manifold, and  therefore $u$ must be constant, since it is superharmonic (see \cite[Corollary 1]{CY}). Hence, $u=0$ on  $M$.  

We also notice that this case is a particular case of Theorem \ref{thm_serrin} below, which then provides an alternative proof.

Now we assume $d \geq 3$.  By the strong maximum principle we immediately have that if $u(x_0)=0$ for some $x_0\in M$ then $u \equiv 0$ in $M$. Hence, we assume that $u>0$ and we define $v$ as in \eqref{v_def}. 

Since $\Ric \geq 0$ and from \eqref{eqPt}, we have
\begin{equation} \label{Alb0}
 \left(t-\frac12\right) P^{t-2} v^{2-n}|\nabla P|^2 + \frac{n-d}{n}P^{t+1}v^{1-n}  \leq \divergence \left( P^{t-1}  v^{2-n}\nabla P \right) 
\end{equation}
in $M$, for any $t \in \R$. Let $\varphi$ be a standard cutoff function as in \eqref{vphi_def}. By multiplying \eqref{Alb0} by $\varphi^\theta$, for some large $\theta$ to be chosen later, and integrating by parts we obtain 
\begin{equation*}
 \left(t-\frac12\right) \int_M \varphi^{\theta} P^{t-2} v^{2-n}|\nabla P|^2 + \frac{n-d}{n} \int_M \varphi^{\theta}  P^{t+1}v^{1-n}  \leq - \theta \int_M  \varphi^{\theta -1} P^{t-1} v^{2-n}\nabla P \cdot \nabla \varphi
\end{equation*}
and hence by Young inequality
\begin{multline*}
 \left(t-\frac12\right) \int_M \varphi^{\theta}  P^{t-2} v^{2-n}|\nabla P|^2 + \frac{n-d}{n} \int_M \varphi^{\theta}  P^{t+1}v^{1-n}  \leq  \int_M (P^{\frac{t-2}{2}} v^{\frac{2-n}{2}} |\nabla P| \varphi^{\frac{\theta}{2}}) (\theta P^{\frac{t}{2}}v^{\frac{2-n}{2}} |\nabla \varphi| \varphi^{\frac{\theta}{2}-1}) \\
 \leq \epsilon \int_M P^{t-2} v^{2-n} |\nabla P|^2 \varphi^\theta  + C(\epsilon) \theta^2 \int_M P^t v^{2-n} |\nabla \varphi |^2 \varphi^{\theta - 2}  \,,
\end{multline*}
i.e.
\begin{equation} \label{Alb1}
 \left(t-\frac12 - \epsilon \right) \int_M \varphi^{\theta}  P^{t-2} v^{2-n}|\nabla P|^2 + \frac{n-d}{n} \int_M \varphi^{\theta}  P^{t+1}v^{1-n}  \leq C(\epsilon) \theta^2 \int_M P^t v^{2-n} |\nabla \varphi |^2 \varphi^{\theta - 2}  \,.
\end{equation}
Now we use Young inequality again to estimate the RHS in \eqref{Alb1}, as follows
\begin{multline*}
C(\epsilon) \theta^2  \int_M P^t v^{2-n} |\nabla \varphi |^2 \varphi^{\theta - 2} = \int_M \left(P^t v^{(1-n)\frac{t}{1+t}} \varphi^{\theta \frac{t}{1+t}}  \right) \left( C(\epsilon) \theta^2  v^{\frac{1-n}{1+t}} v |\nabla \varphi |^2 \varphi^{\frac{\theta}{1+t} - 2} \right) \\
\leq \epsilon   \int_M P^{t+1} v^{1-n} \varphi^{\theta} + \tilde C(\epsilon, \theta, t) \int_M v^{2-n+t} |\nabla \varphi |^{2(t+1)} \varphi^{\theta - 2(1+t)} 
\end{multline*}
and from \eqref{Alb1} we obtain
\begin{multline} \label{Alb2}
 \left(t-\frac12 - \epsilon \right) \int_M \varphi^{\theta}  P^{t-2} v^{2-n}|\nabla P|^2 + \left( \frac{n-d}{n} - \epsilon \right) \int_M \varphi^{\theta}  P^{t+1}v^{1-n} \\
  \leq \tilde C(\epsilon, \theta, t) \int_M v^{-(n-2-t)} |\nabla \varphi |^{2(t+1)} \,.
\end{multline}
We choose 
$$
t=n-2 \quad  \textmd{ and } \quad 0< \epsilon <  \frac{n-d}{n} \,, 
$$
so that
$$
t-\frac12 - \epsilon   > 0  \quad \textmd{ and } \quad  \frac{n-d}{n} - \epsilon > 0 
$$
for any $3 \leq d < n$. By choosing $\theta \geq 2(t+1)$ and since $\varphi$ satisfies \eqref{vphi_def}, then \eqref{Alb2} implies
\begin{equation*} 
 \left(t-\frac12 - \epsilon \right) \int_M \varphi^{\theta}  P^{t-2} v^{2-n}|\nabla P|^2 + \left( \frac{n-d}{n} - \epsilon \right) \int_M \varphi^{\theta}  P^{t+1}v^{1-n} 
  \leq  C_*(\epsilon, \theta, t) R^{-(n-2)-(n-d)} \,.
\end{equation*}
In particular, this implies that
$$
\int_{B_R} P^{t+1}v^{1-n} \leq  C R^{-(n-2)-(n-d)}  \,,
$$
where $C$ does not depend on $R$. By letting $R \to +\infty$, we obtain that $P^{t+1}v^{1-n} \equiv 0$ in $M$, which is a contradiction. Hence, $u$ must vanish everywhere in $M$.
%\begin{equation*} 
% \left(t-\frac12 - \epsilon \right) \int_M \varphi^{\theta}  P^{t-2} v^{2-n}|\nabla P|^2 + \left( \frac{n-d}{n} - \epsilon \right) \int_M \phi^{\theta}  P^{t+1}v^{1-n} 
%  \leq  C_*(\epsilon, \theta, t) R^{-2(t+1)}  \int_M v^{-(n-2-t)}
%\end{equation*}
%and \eqref{(ii)} with $q= n-2-t$ gives  
%\begin{multline*} 
% \left(t-\frac12 - \epsilon \right) \int_{B_R}   P^{t-2} v^{2-n}|\nabla P|^2 + \left( \frac{n-d}{n} - \epsilon \right) \int_{B_R} P^{t+1}v^{1-n} \\
%  \leq C_*(\epsilon, \theta, t) R^{-2(t+1)+d-(n-2-t)} = C_*(\epsilon, \theta, t) R^{d-n-t} \,.
%\end{multline*}
%Since $d-n-t<0$, by letting $R\to +\infty$ we obtain that $ P^{t+1}v^{1-n} \equiv 0$ in $M$, which is a contradiction.
\end{proof}

Now we consider the case of Theorem \ref{thm_subcritical_Ricci_neg}, where we assume that $d \geq 3$ and that the manifold enjoys a Sobolev type inequality \eqref{Sobolev_nu}, i.e.
\begin{equation*} 
\|\nabla \psi \|_{L^2(M)} \geq S \|\psi\|_{L^{\nu}(M)} \qquad \forall \,  \psi \in C_c^1(M),
\end{equation*}
for some $\nu >2$. We recall that in this case we do not restrict to manifolds with nonnegative Ricci curvature, but we assume that the negative part of the Ricci tensor satisfies \eqref{Ricci_bound_neg}, i.e.
\begin{equation*} 
\|\Ric_{-}\|_{L^{\frac{\nu}{\nu-2}}(M)}  \leq  \frac{n-d}{3(n-2)^2} S^2  \,,
\end{equation*}
(see Remark \ref{remark_intro}). Moreover, we assume the volume growth condition \eqref{Vol_bound_in_p} which is given by
\begin{equation}  \label{vg_cond}
V_g(B_R) = O( R^{2(n-1)})
\end{equation}
as $R\to \infty$ (again, see Remark \ref{remark_intro}).

\begin{proof}[Proof of Theorem \ref{thm_subcritical_Ricci_neg}]
By the strong maximum principle we have that if $u(p)=0$ for some $p\in M$ then $u \equiv 0$ in $M$. 
Hence, by contradiction, we assume that $u>0$ in $M$ and we define $v$ as in \eqref{v_def}. 

Let $\vp$ be a standard cutoff function as in \eqref{vphi_def} and let $\theta >0$ to be chosen later. From  \eqref{eq_v} and \eqref{eqPt} we get
\begin{multline*} 
\int_M \vp^\theta \divergence \left( P^{t-1}  v^{2-n}\nabla P \right) \\ \geq \int_M \left(t-\frac12\right) P^{t-2} v^{2-n}|\nabla P|^2 \vp^\theta + n \int_M P^{t-1}v^{1-n}  \left\{ \Ric(\nabla v,\nabla v) + \frac{n-d}{n^2}  P^2 \right\}  \vp^\theta 
\end{multline*}
for any $t \in \R$ and, by using the divergence theorem, we find
\begin{multline*} 
 \int_M \left(t-\frac12\right) P^{t-2} v^{2-n}|\nabla P|^2 \vp^\theta   + \frac{n-d}{n} \int_M P^{t+1}v^{1-n}  \vp^\theta  \\ \leq  n \int_M \Ric_{-} \vp^\theta P^{t-1}v^{1-n}  |\nabla v|^2 - \int_M \theta \vp^{\theta-1}  P^{t-1}  v^{2-n}\nabla P \cdot \nabla \vp \\
\leq 2 \int_M \Ric_{-} \vp^\theta P^{t}v^{2-n}  - \int_M \theta \vp^{\theta-1}  P^{t-1}  v^{2-n}\nabla P \cdot \nabla \vp  \,.
\end{multline*}
where in the last inequality we used that  
\begin{equation} \label{nabla v leq}
|\nabla v|^2 < \frac{2}{n} v P \,,
\end{equation}
which follows from $vP = \frac{n}{2}|\nabla v|^2 + \frac{2}{n-2}$.
From H\"older inequality with exponents $(\frac{\nu}{2},\frac{\nu}{\nu-2})$ we have that
\begin{equation*}
2 \int_M \Ric_{-} \vp^\theta  P^{t}v^{2-n} \leq 2 \|\Ric_{-}\|_{\frac{\nu}{\nu-2}}  \| v^{2-n}\vp^\theta  P^{t} \|_{\frac{\nu}{2}} \,,
\end{equation*}
and we obtain
\begin{multline*} %\label{start1} 
\left(t-\frac12\right) \int_M  P^{t-2} v^{2-n}|\nabla P|^2 \vp^\theta   + \frac{n-d}{n} \int_M P^{t+1}v^{1-n}  \vp^\theta  \\
 \leq 2 \|\Ric_{-}\|_{\frac{\nu}{\nu-2}}  \| v^{2-n}\vp^\theta  P^{t} \|_{\frac{\nu}{2}}  - \int_M \theta \vp^{\theta-1}  P^{t-1}  v^{2-n}\nabla P \cdot \nabla \vp  \,.
\end{multline*}
Hence, Young inequality yields
\begin{multline} \label{start2} 
\left(t-\frac12 - \epsilon \right) \int_M  P^{t-2} v^{2-n}|\nabla P|^2 \vp^\theta   + \frac{n-d}{n} \int_M P^{t+1}v^{1-n}  \vp^\theta  \\
 \leq 2 \|\Ric_{-}\|_{\frac{\nu}{\nu-2}}  \| v^{2-n}\vp^\theta  P^{t} \|_{\frac{\nu}{2}} +  \frac{\theta^2}{4\epsilon} \int_M P^t \vp^{\theta-2} v^{2-n} |\nabla \vp|^2  \,.
\end{multline}
Since 
$$
 P^{t-2} v^{2-n}|\nabla P|^2 \vp^\theta = |v^{1-\frac{n}{2}} \vp^{\frac{\theta}{2}} P^{\frac{t}{2}-1}\nabla P|^2 =\frac{4}{t^2} |v^{1-\frac{n}{2}} \vp^{\frac{\theta}{2}} \nabla P^{\frac{t}{2}}|^2 =\frac{4}{t^2} v^{2-n} \vp^{\theta} |\nabla P^{\frac{t}{2}}|^2 
$$
and 
$$
 |\nabla (v^{1-\frac{n}{2}} \vp^{\frac{\theta}{2}}  P^{\frac{t}{2}})|^2 \leq 3 \left(1-\frac{n}{2} \right)^2 v^{-n} \vp^\theta P^t |\nabla v|^2+ 3 \frac{\theta^2}{4} \vp^{\theta-2} v^{2-n} P^t |\nabla \vp|^2 + 3\ v^{2-n} \vp^{\theta} |\nabla P^{\frac{t}{2}}|^2 \,,
$$
i.e.,
$$
\frac{4}{t^2} v^{2-n} \vp^{\theta} |\nabla P^{\frac{t}{2}}|^2  \geq \frac{4}{3t^2} |\nabla (v^{1-\frac{n}{2}} \vp^{\frac{\theta}{2}}  P^{\frac{t}{2}})|^2  -  \frac{4}{t^2} \left(1-\frac{n}{2} \right)^2 v^{-n} \vp^\theta P^t |\nabla v|^2 - \frac{\theta^2}{t^2} \vp^{\theta-2} v^{2-n} P^t |\nabla \vp|^2 \,,
$$
we find 
\begin{multline*}
 \left(t-\frac12 -\epsilon  \right) P^{t-2} v^{2-n}|\nabla P|^2 \vp^\theta \geq  \left(t-\frac12 -\epsilon \right)\frac{4}{3t^2} |\nabla (v^{1-\frac{n}{2}} \vp^{\frac{\theta}{2}}  P^{\frac{t}{2}})|^2 \\
 -  \left(t-\frac12 -\epsilon \right) \frac{4}{t^2} \left(1-\frac{n}{2} \right)^2 v^{-n} \vp^\theta P^t |\nabla v|^2   -  \left(t-\frac12 -\epsilon \right) \frac{\theta^2}{t^2} \vp^{\theta-2} v^{2-n} P^t |\nabla \vp|^2 \,,
\end{multline*}
provided $t> \frac12$. From \eqref{start2} we obtain  
\begin{multline*} 
\left(t-\frac12 -\epsilon \right)\frac{4}{3t^2} \int_M   |\nabla (v^{1-\frac{n}{2}} \vp^{\frac{\theta}{2}}  P^{\frac{t}{2}})|^2 + \frac{n-d}{n} \int P^{t+1}v^{1-n}     \vp^\theta \\  
\leq  2 \|\Ric_{-}\|_{\frac{\nu}{\nu-2}}  \| v^{2-n}\vp^\theta  P^{t} \|_{\frac{\nu}{2}} + \frac{\theta^2}{4\epsilon} \int_M P^t \vp^{\theta-2} v^{2-n} |\nabla \vp|^2 \\
+ \left(t-\frac12 -\epsilon \right) \frac{4}{t^2} \left(1-\frac{n}{2} \right)^2 \int_M v^{-n} \vp^\theta P^t |\nabla v|^2  + \left( t-\frac12-\epsilon \right) \frac{\theta^2}{t^2}  \int_M P^t \vp^{\theta-2} v^{2-n} |\nabla \vp|^2   
\end{multline*}
and from \eqref{nabla v leq} we find
\begin{multline}  \label{restart}
\left(t-\frac12 -\epsilon \right)\frac{4}{3t^2} \int_M   |\nabla (v^{1-\frac{n}{2}} \vp^{\frac{\theta}{2}}  P^{\frac{t}{2}})|^2 +\frac{1}{n} \left[ n-d -  \left(t-\frac12 -\epsilon \right) \frac{2(n-2)^2}{t^2} \right]  \int_M P^{t+1}v^{1-n}     \vp^\theta \\  
\leq  2 \|\Ric_{-}\|_{\frac{\nu}{\nu-2}}  \| v^{2-n}\vp^\theta  P^{t} \|_{\frac{\nu}{2}}  +   \left[ \frac{1}{4\epsilon} +  \left(t-\frac12-\epsilon \right) \frac{1}{t^2} \right]\theta^2 \int_M P^t \vp^{\theta-2} v^{2-n} |\nabla \vp|^2  \,.
\end{multline}
Sobolev inequality \eqref{Sobolev_nu} yields
\begin{multline}  \label{restart}
2\left[\left(t-\frac12 -\epsilon \right)\frac{2}{3t^2}S^2  -  \|\Ric_{-}\|_{\frac{\nu}{\nu-2}} \right] \| v^{2-n}\vp^\theta  P^{t} \|_{\frac{\nu}{2}}  \\
+\frac{1}{n} \left[ n-d -  \left(t-\frac12 -\epsilon \right) \frac{2}{t^2}(n-2)^2 \right]  \int_M P^{t+1}v^{1-n}     \vp^\theta \\  
\leq     \left[ \frac{1}{4\epsilon} +  \left(t-\frac12-\epsilon \right) \frac{1}{t^2} \right]\theta^2 \int_M P^t \vp^{\theta-2} v^{2-n} |\nabla \vp|^2  \,.
\end{multline}
Now we choose the parameters in \eqref{restart} in such a way that 
\begin{equation} \label{conz1}
\begin{cases}
\left(t-\frac12 - \epsilon \right)\frac{2}{3t^2} S^2  -  \|\Ric_{-}\|_{\frac{\nu}{\nu-2}}  \geq 0 & \\
n-d-  \left(t-\frac12 - \epsilon \right) \frac{2}{t^2} \left(n-2\right)^2 >0  \,.& 
\end{cases}
\end{equation}
In particular, we first choose $t$ and $\epsilon$ such that 
\begin{equation} \label{conz2}
\begin{cases}
t-\frac12 - \epsilon > 0 & \\
n-d  -  \left(t-\frac12 - \epsilon \right) \frac{2}{t^2} \left( n-2 \right)^2 > 0 \,. &
\end{cases}
\end{equation}
In particular, we set
\begin{equation} \label{t_ep_def}
t= n-2 \quad \textmd{ and } \quad \epsilon= n-2- \frac12  - \frac{n-d}{4} \,,
\end{equation}
and we notice that, since $n>d \geq 3,$ then $t>1$ and $\epsilon > \frac12$. Moreover, such a choice of $t$ and $\epsilon$ implies that 
$$
t -\epsilon - \frac12 = \frac{n-d}{4} > 0 
$$
and 
$$
n-d  -  \left(t-\frac12 - \epsilon \right) \frac{2}{t^2} \left( n-2 \right)^2 = \frac{n-d}{2} > 0
$$
and \eqref{conz2} is satisfied. 

Moreover, by noticing that
$$ 
\left(t-\frac12 - \epsilon \right)\frac{2}{3t^2} S^2  = \frac{n-d}{3(n-2)^2}   S^2 \,,
$$
we also have that the first inequality in \eqref{conz1} is satisfied, thanks to the assumption \eqref{K-cond} on $\|\Ric_{-}\|_{\frac{\nu}{\nu-2}} $. 

Now that \eqref{conz1} is fulfilled, from \eqref{restart} we find 
\begin{equation} \label{ft5} 
c_1 \| v^{2-n}\vp^\theta  P^{t} \|_{\frac{\nu}{2}}  + c_2 \int_M P^{t+1}v^{1-n}  \vp^\theta   
\leq    \theta^2  \int_M P^t \vp^{\theta-2} v^{2-n} |\nabla \vp|^2  \,,
\end{equation}
where $c_1=c_1(n,d,S,\|\Ric_{-}\|_{\frac{\nu}{\nu-2}}) \geq 0$ and $c_2=c_2(n,d)>0$. Now we estimate the last term in \eqref{ft5}. We write
\begin{equation} \label{preHolder}
\theta^2 \int_M P^t \vp^{\theta-2} v^{2-n} |\nabla \vp|^2 = \int_M (  P^t  v^{(1-n)\frac{t}{1+t}} \vp^{\theta \frac{t}{1+t}}   ) ( \theta^2 v^{\frac{1-n}{1+t}} v  \vp^{\frac{\theta}{1+t}-2} |\nabla \vp|^2 )
\end{equation}
and, by using Young's inequality with exponents $(\frac{t+1}{t},t+1)$, we find 
$$
\theta^2 \int_M P^t \vp^{\theta-2} v^{2-n} |\nabla \vp|^2  \leq  \frac{c_2}{2} \int_M  P^{t+1}  v^{(1-n)} \vp^{\theta } + C(c_2,\theta) \int_M  v^{2-n + t}   \vp^{\theta-2(1+t)} |\nabla \vp|^{2(t+1)} \,.
$$
By choosing $\theta > 2(1+t)$ and using the properties of $\vp$, from \eqref{ft5} we find
\begin{equation} \label{ft51} 
c_1 \| v^{2-n}\vp^\theta  P^{t} \|_{\frac{\nu}{2}}  + \frac{c_2}{2} \int_M P^{t+1}v^{1-n}  \vp^\theta 
\leq   C(c_2,\theta) R^{-2(t+1)} \int_{B_{2R} \setminus B_R}  v^{2-n + t} \leq  C(c_2,\theta) R^{-2(n-1)} V_g (B_{2R})   \,,
\end{equation}
where in the last inequality we used that $t=n-2$.  Now, it remains to use the volume growth condition \eqref{vg_cond} to conclude. 

If we assume that $V_g(B_{2R}) = o(R^{-2(n-1)})$ then, by letting $R \to \infty$, from \eqref{ft51} we immediately obtain that 
\begin{equation*} 
0 \leq c_1 \| v^{2-n}   P^{t} \|_{\frac{\nu}{2}}  + \frac{c_2}{2} \int_M P^{t+1}v^{1-n} \leq 0 
\end{equation*}
so that $P=0$, which gives a contradiction and we conclude.

If $V_g(B_{2R}) = O(R^{-2(n-1)})$ then \eqref{ft51} yields 
\begin{equation*} 
c_1 \| v^{2-n}   P^{t} \|_{\frac{\nu}{2}}  + \frac{c_2}{2} \int_M P^{t+1}v^{1-n} < +\infty.
\end{equation*}
which implies that 
\begin{equation} \label{int_bound} 
\lim_{R\to \infty} \int_{B_{2R}\setminus B_R} P^{t+1}v^{1-n} =0 \,.
\end{equation}
Since $\supp \nabla \vp \subset B_{2R} \setminus B_R$ and from \eqref{ft5} we find
$$ 
c_2 \int_M P^{t+1}v^{1-n}  \vp^\theta   \leq    \theta^2  \int_M P^t \vp^{\theta-2} v^{2-n} |\nabla \vp|^2 =  \theta^2  \int_{B_{2R} \setminus B_R} P^t \vp^{\theta-2} v^{2-n} |\nabla \vp|^2
$$
and, by arguing as in \eqref{preHolder} and using H\"older inequality, we obtain
\begin{equation} \label{quasi}
c_2 \int_M P^{t+1}v^{1-n}  \vp^\theta   \leq    \theta^2  \left( \int_{B_{2R} \setminus B_R}  P^{t+1}  v^{(1-n)} \vp^{\theta } \right)^{\frac{t}{t+1}} \left( \int_{B_{2R} \setminus B_R}  v^{2-n + t}   \vp^{\theta-2(1+t)} |\nabla \vp|^{2(t+1)} \right)^{\frac{1}{t+1}} \,.
\end{equation}
By repeating the argument done in \eqref{ft51}, from \eqref{quasi} we obtain
\begin{equation*}  
c_2 \int_M P^{t+1}v^{1-n}  \vp^\theta   \leq    \theta^2  \left( \int_{B_{2R} \setminus B_R}  P^{t+1}  v^{(1-n)} \vp^{\theta } \right)^{\frac{t}{t+1}} \left(C V_g(B_{2R}) R^{-2(n-1)} \right)^{\frac{1}{t+1}} \,.
\end{equation*}
Since $V_g(B_{4R}) = O(R^{-2(n-1)})$ then  we find that
\begin{equation*}  
\int_{B_R} P^{t+1}v^{1-n} \leq \int_M P^{t+1}v^{1-n}  \vp^\theta   \leq   C  \left( \int_{B_{2R} \setminus B_R}  P^{t+1}  v^{(1-n)} \right)^{\frac{t}{t+1}}  
\end{equation*}
and, by taking the limit as $R\to \infty$ and using \eqref{int_bound}, we find that
\begin{equation*}  
\int_{M} P^{t+1}v^{1-n}  =  0  
\end{equation*}
which gives the desired contradiction.
\end{proof}

\section{Integral estimates for solutions to  \eqref{eq_v}} \label{section_integral_estimates}

In this section we prove some integral estimates for solutions to  \eqref{eq_v}, which will be used for proving Theorems \ref{thm_critical} and \ref{thm_Liouville}. 

\begin{lemma}\label{Lemma1}
Let $(M,g)$ be a Riemannian manifold of dimension $d \geq 2$. Let $v$ be a positive solution to \eqref{eq_v} and let $n$ be given by \eqref{n_def}. Then, for any $q\in\R$ and $\psi\in  C^{0,1}_c(M)$, it holds
\begin{equation} \label{eq_Lemma1}
(\frac{n}{2}+1-q)\int_{M} v^{-q}|\nabla v|^2\psi+ c_n\int_{M} v^{-q}\psi = -\int_{M} v^{1-q}\nabla v \cdot \nabla \psi \,,
\end{equation}
where $c_n$ is given by \eqref{cn_def}.
\end{lemma}

\begin{proof}
From \eqref{eq_v} and \eqref{P_def12} and by integrating by parts, we get
\begin{equation*}
\int_{M} v^{1-q}\psi P =\int_{M}   v^{1-q}\psi \Delta v=(q-1)\int_{M} v^{-q}|\nabla v|^2\psi -\int_{M} v^{1-q}\nabla v\cdot \nabla \psi \,.
\end{equation*}
Using again  \eqref{P_def12} we immediately get \eqref{eq_Lemma1}.
%$$
%(\frac{n}{2}+1-q)\int_{M} v^{-q}|\nabla v|^2\psi+\frac{2}{n-2}\int_{M} v^{-q}\psi = -\int_{M} v^{1-q}\nabla v \cdot \nabla \psi \,.
%$$
\end{proof}

\begin{corollary}\label{Corollary}
Let $(M,g)$ be a Riemannian manifold of dimension $d \geq 2$. Let $v$ be a positive solution to \eqref{eq_v} and let $n$ be given by \eqref{n_def}. Let $q,m\in\R$, $\psi\in  C^{0,1}_c(M)$ and let $P$ be given by \eqref{P_def12}. Then it holds
\begin{multline} \label{eq_Corollary} 
(\frac{n}{2}+1-q)\int_{M} v^{-q}P^m|\nabla v|^2\psi+ c_n\int_{M} v^{-q}P^m\psi  \\
= -\int_{M} v^{1-q}P^m\nabla v \cdot \nabla \psi - m \int_{M} v^{1-q}P^{m-1}\nabla v \cdot \nabla P \psi \,,
\end{multline}
where $c_n$ is given by \eqref{cn_def}.
\end{corollary}

\begin{proof}
We recall that, thanks to the regularity properties of $v$, $P$ is smooth and we can choose $\psi P^m$ in place of $\psi$ in Lemma \ref{Lemma1}, which immediately gives \eqref{eq_Corollary}.
\end{proof}

In the following lemma we provide integral estimates which play a crucial role in the proof of Theorem \ref{thm_critical}. These are rather standard tools to study positive solutions of $p$-Laplace type equations involving power-like nonlinearities (see for instance \cite{GS, SerrinZou, CatinoMonticelliRoncoroni, OU, Vetois}).

\begin{lemma}\label{Lemma2bis}
Let $(M,g)$ be a complete, connected  Riemannian manifold of dimension $d \geq 2$ and assume that $n>d$ if $d=2$, where $n$ is given by \eqref{n_def}. Let $v$ be a positive solution to \eqref{eq_v}. For any $R>0$ we have 
\begin{itemize}
\item[(i)] if $2\le q<\frac{n}{2} + 1$ then 
\begin{equation} \label{(i)}
\int_{B_R} v^{-q}|\nabla v|^2+\int_{B_R} v^{-q}\le C \, V_g(B_{2R})R^{-q}  \,.
\end{equation}
%\todo{da togliere: where, for a fixed point $p \in M$,  $V_g(B_r)$ denotes the volume of the geodesic ball $B_r = B_r(p)$ of center $p$ and radius $r>0$. }

\item[(ii)] If $ 0 \leq q\le \frac{n}{2}+1$ then
\begin{equation} \label{(ii)}
\int_{B_{R}} v^{-q}\le  C \, V_g(B_{2R})R^{-q} \,.
\end{equation}
%\todo{da togliere: where, for a fixed point $p \in M$,  $V_g(B_r)$ denotes the volume of the geodesic ball $B_r = B_r(p)$ of center $p$ and radius $r>0$.}
\end{itemize}
In both cases, $V_g(B_r)$ is the volume of the geodesic ball $B_r = B_r(x_0)$ centered at a point $x_0 \in M$ and of radius $r>0$, and $C=C(n,q)$ is a positive constant depending only on $n$ and $q$.    
\end{lemma}

\begin{proof}

Let $\varphi$ be a cutoff function as in \eqref{vphi_def}.

\noindent (i) Assume $\theta >q$. By letting $\psi = \varphi^\theta$ in Lemma \ref{Lemma1} and using Cauchy-Schwarz inequality, we have
\begin{equation} \label{sotto}
(\frac{n}{2}+1-q)\int_{M} v^{-q}|\nabla v|^2\varphi^\theta+ c_n \int_{M} v^{-q}\varphi^\theta \leq  \theta\int_{M} \varphi^{\theta-1}v^{1-q}|\nabla v| |\nabla \varphi| \,.
\end{equation}
 Let $\epsilon >0$ be a (small) constant to be chosen later. By using Young's inequality we have
\begin{equation} \label{yyoung}
\varphi^{\theta-1} v^{1-q}|\nabla v| |\nabla\varphi| \leq \epsilon\varphi^\theta v^{-q}|\nabla v|^2+\frac{1}{\epsilon}  |\nabla\varphi|^2 \varphi^{\theta-2}v^{2-q} \,,
\end{equation}
and a further application of Young's inequality, this time with exponents $\big(\frac{q}{q-2},\frac{q}{2}\big)$ and assuming $q>2$, yields
\begin{equation} \label{young2volte}
\varphi^{\theta-1} v^{1-q}|\nabla v| |\nabla\varphi| \leq  \epsilon\varphi^\theta v^{-q}|\nabla v|^2 + \epsilon \varphi^\theta v^{-q} + C(\epsilon,q)|\nabla \varphi|^q \varphi^{\theta -q} \, ,
\end{equation} 
where $C=C(\epsilon,q)$ is a positive constant depending only on $ \epsilon$ and $q$.   
We notice that if $q=2$ then \eqref{young2volte} immediately follows from \eqref{yyoung}.

From \eqref{sotto}, \eqref{young2volte} and by using \eqref{vphi_def}, we get
\begin{multline}\label{comp30}
(\frac{n}{2}+1-q)\int_{M} v^{-q}|\nabla v|^2\varphi^\theta+ c_n \int_{M} v^{-q}\varphi^\theta \\ \le  
\theta\epsilon\int_{B_{2R} \setminus B_R} \varphi^\theta v^{-q}|\nabla v|^2+ \theta \epsilon\int_{B_{2R} \setminus B_R} v^{-q} \varphi^{\theta} + \theta C(\epsilon,q) \int_{B_{2R} \setminus B_R}  |\nabla \varphi|^q \varphi^{\theta-q} \,.
\end{multline}
i.e.
\begin{equation}\label{comp30bis}
(\frac{n}{2}+1-q-\theta \epsilon)\int_{M} v^{-q}|\nabla v|^2\varphi^\theta+ (c_n -\theta \epsilon) 
\int_{M} v^{-q}\varphi^\theta \le  \theta C(\epsilon,q) R^{-q} \int_{B_{2R} \setminus B_R}  \varphi^{\theta-q} \,.
\end{equation}
Since $\theta >q$ and $(\frac{n}{2}+1-q)>0$ and by taking $\epsilon>0$ small enough, we obtain
$$
\int_{B_R}  v^{-q}|\nabla v|^2+\int_{B_R} v^{-q}\le C \, V_g(B_{2R})R^{-q},
$$
where $C=C(n,q)$ is a positive constant depending only on $n$ and $q$.   

\bigskip

Now we prove (ii). We first notice that the case $q=0$ is trivially true, while if $2 \leq q<n/2+1$ then \eqref{(ii)} follows from \eqref{(i)}. 

In the case $0< q < 2$, \eqref{(ii)} follows by a direct application of H\"older inequality and part (i) of the assertion. Indeed, let $0< s < 2$. We apply H\"older inequality with exponents $(\frac{q}{s},\frac{q}{q-s})$, with $2< q<n/2+1$, and get
\begin{multline}
\int_{B_{R}} v^{-s}\le\left(\int_{B_{R}}v^{-q}\right)^{s/q} V_g(B_{R})^{(q-s)/q} %\\ 
\le (C \, V_g(B_{2R})R^{-q})^{\frac{s}{q}} V_g(B_{R})^{(q-s)/q} \leq  C ^{\frac{s}{q}}  V_g(B_{2R})R^{-s}, 
\end{multline}
where in the second  inequality we used \eqref{(i)} for $2< q<n/2+1$. Therefore, the estimate \eqref{(ii)} holds for $0\le q<n/2+1$. 

It remains to prove \eqref{(ii)} for $q=n/2+1$. 
%Notice that we can assume that $n>2$, since for $n=2$ we would have $q=2$ and the assertion follows from \eqref{(i)}.  
From Lemma \ref{Lemma1} we have 
\begin{equation} \label{int_qn2}
c_n \int_{M} v^{-\frac{n}{2}-1}\varphi \le \frac{1}{R} \int_{B_{2R}\setminus B_R} v^{-\frac{n}{2}}|\nabla v| \,.
\end{equation}
Now we estimate the right hand side of the last inequality. By writing $v^{-\frac{n}{2}}= v^{-\frac{n}{4}-\frac12+\frac{\epsilon}{2}}v^{-\frac{n}{4}+\frac12-\frac{\epsilon}{2}}$ and using Cauchy-Schwarz inequality we find
$$
\int_{B_{R}} v^{-\frac{n}{2}}|\nabla v| \le \left(\int_{B_{R}} v^{\epsilon-\frac{n}{2}-1}|\nabla v|^2\right)^{\frac12} \left(\int_{B_{R} }v^{-\frac{n}{2}+1-\epsilon}\right)^{\frac12}.
$$
By choosing $0<\epsilon<\frac{n-2}{2}$ on the r.h.s.  we can apply \eqref{(i)} for $2< q<\frac{n}{2}+1$ on the first term, and \eqref{(ii)} for $0<q<\frac{n}{2}+1$ on the second term to get
\beq\label{comp10}
\int_{B_{R}} v^{-\frac{n}{2}}|\nabla v| \le  \left( C  V_g(B_{2R})R^{-(-\epsilon+\frac{n}{2}+1)} \right) ^{\frac{1}{2}}
 \left( C  V_g(B_{2R})R^{-(\epsilon+\frac{n}{2}-1)} \right) ^{\frac{1}{2}}.
\eeq
Inequality \eqref{(ii)} then follows from \eqref{int_qn2} and \eqref{comp10}.
\end{proof}

\begin{remark} \label{remark_asymp}
It is clear that from \eqref{P_def12} and \eqref{(i)} we obtain that 
\begin{equation} \label{(i)bis}
\int_{B_R} v^{\gamma} P  \le C \, V_g(B_{2R})R^{\gamma-1} 
\end{equation}
for any $ -\frac{n}{2} < \gamma \leq -1$. These estimate will be used later.
\end{remark}

We notice that part (ii) of Lemma \ref{Lemma2bis} implies the following straightforward classification result, where only a volume growth is assumed on the manifold $M$.

\begin{theorem} \label{thm_serrin}
Let $(M,g)$ be a complete, connected Riemannian manifold of dimension $d \geq 2$ and let $u$ be a nonnegative solution to \eqref{eq_semilinear1}, where $1<p \leq p_S$ when $d\geq 3$ and $1<p<+\infty$ for $d=2$.

If
\begin{equation} \label{hyp-vol_Serrin}
V_g(B_R) = {o(R^{\frac{2p}{p-1}})}
%o(R^{\frac{2p}{p-2}})
\end{equation}
then $u=0$. 

In particular, if 
\begin{equation} \label{vol_gamma}
V_g(B_R) \leq C R^\gamma \,,
\end{equation}
for some $C,\gamma>0$ then $u=0$ whenever $p<p_M$, where
\begin{equation} \label{p_serrin}
p_M= 
\begin{cases}
\dfrac{\gamma}{(\gamma-2)^+}  & \textmd{ if } d \geq 3 \\
+\infty & \textmd{ if } d=2 \,.
\end{cases}
\end{equation}
\end{theorem}

\noindent 
We notice that, when $\gamma=d \geq 3$ in \eqref{vol_gamma} then $p_M=d/(d-2)^+,$ which is the well-known \emph{Serrin's exponent}.

\begin{proof}[Proof of Theorem \ref{thm_serrin}]
By the strong maximum principle either $u = 0$ on $M$ or $ u>0$ on $M$. Let us prove that under the assumption \eqref{hyp-vol_Serrin} the latter case is not possible. For contradiction,  assume that $u>0$ in $M$ and define $v$ as in \eqref{v_def}, then $\liminf_{R \to \infty} \frac{V_g(B_{R})}{R^{\frac{n}{2}+1}} = \liminf_{R \to \infty} \frac{V_g(B_{R})}
{R^{\frac{2p}{p-1}}} >0$, thanks to part (ii) of Lemma \ref{Lemma2bis}. The latter contradicts \eqref{hyp-vol_Serrin}. Hence $u=0$ on $M$.

\end{proof}

\section{The critical Laplace equation for manifolds with nonnegative Ricci} \label{section5}
In this section we assume that $\Ric \geq 0$ and we prove Theorem \ref{thm_critical}. We notice that throughout this section, we consider solutions to \eqref{eq_semilinear1} and hence in this case $n=d$ (see \eqref{n_def}). However, we prefer to use both and $n$ and $d$, since $n$ appears in quantities which are related to the equation while the dimension $d$ plays a role only in the asymptotic volume ratio. This will be useful for future works and generalizations.  

We firstly give two preliminary propositions.

\begin{proposition}\label{mainineq}
Let $(M,g)$ be a complete, connected  Riemannian manifold of dimension $d\geq 3$ with $\Ric \geq 0$. Let $n$ and $P$ be given by  \eqref{n_def} and \eqref{P_def12}, respectively, and assume that $n=d$. Let $\varphi$ be a standard cutoff function as in \eqref{vphi_def}. If $v$ is a positive solution to \eqref{eq_v} then 
\beq\label{comp26}
\frac{1}{\theta^2}\left(t-\frac12\right)^2 \int_{M}  P^{t-2}v^{2-n}|\nabla P|^2\vp^\theta  \leq 
\frac{1}{R^2} \int_{B_{2R} \setminus B_R}\vp^{\theta-2}v^{2-n}P^t  \,,
\eeq
for any $t>\frac{1}{2}$ and $\theta\ge 2$.
\end{proposition}

\begin{proof}
We first notice that, since $n= d$ and $\Ric \geq 0$, from \eqref{eqPt} we have
\begin{equation*}  
\vp^\theta \divergence \left( P^{t-1}  v^{2-n}\nabla P \right)  \geq \left(t-\frac12\right) P^{t-2} v^{2-n}|\nabla P|^2 
\vp^\theta  \,,    
\end{equation*}
and an integration by parts yields
\begin{multline} \label{intbyparts1}
\left(t-\frac12\right) \int_{M}  P^{t-2}v^{2-n}|\nabla P|^2\vp^\theta  \le - \theta\int_{M} \vp^{\theta-1}  v^{2-n}P^{t-1}\nabla P \cdot \nabla\vp\\
\le \theta \left(\int_{M} \vp^\theta v^{2-n}P^{t-2}|\nabla P|^2\right)^\frac{1}{2}\left(\int_{M} \vp^{\theta-2}v^{2-n}P^t |\nabla \vp|^2\right)^\frac{1}{2}, 
\end{multline}
where in the last inequality we used H\"older's inequality. Hence,  inequality \eqref{comp26} follows since $t > \frac{1}{2}$ and $|\nabla \vp | \leq \frac{1}{R}$ in $B_{2R} \setminus B_R$. 
%\begin{equation*}
%\left(t-\frac12\right)^2 \int_{M}  P^{t-2}v^{2-n}|\nabla P|^2\varphi^\theta  \leq \theta^2 \int_{B_R \setminus B_r}
%\varphi^{\theta-2}v^{2-n}P^t|\nabla \varphi|^2  \,,
%\end{equation*}
%which gives \eqref{comp26}.   
\end{proof}

\begin{proposition}\label{prop_ngrande}
Let $(M,g)$ be a complete, connected  Riemannian manifold of dimension $d\geq 3$ with $\Ric \geq 0$. Let $n$ and $P$ be given by  \eqref{n_def} and \eqref{P_def12}, respectively, and assume that $n=d$. Let $\varphi$ be a standard cutoff function as in \eqref{vphi_def}. If $v$ is a positive solution to \eqref{eq_v}, then for any $q<\frac{n}{2} +1$, $\gamma>1$, $ \theta \geq \max(2,\gamma)$ and $t,m\in \R$ there exists $C>0$ such that
\begin{multline}\label{vac1}
\int_{B_{2R}} v^{1-q} P^{m+1} \varphi^{ \theta} \leq \frac{C}{R^\gamma} \int_{B_{2R} \setminus B_R} \varphi^{ \theta} v^{1-q+ \frac{\gamma}{2}} P^{m+1-\frac{\gamma}{2}} \\ 
+ \frac{C}{R^2}\int_{B_{2R}} v^{n-2q+1} P^{2m-t+1} \varphi^{\theta-2} + \epsilon R^2 \int_{B_{2R}} P^{t-2} v^{2-n} |\nabla P|^2 \varphi^{\theta +2} \,,
\end{multline}
for any positive $\epsilon=\epsilon(n,q)$ sufficiently small. 
\end{proposition}

\begin{proof}
From Corollary \ref{Corollary}, by choosing $\psi = \varphi^\theta$ we know that
\begin{multline*} 
(\frac{n}{2}+1-q)\int_{M} v^{-q}P^m|\nabla v|^2\varphi^\theta+ \frac{2}{n-2} \int_{M} v^{-q}P^m\varphi^\theta  \\
= -\theta \int_{M} \varphi^{\theta-1} v^{1-q}P^m\nabla v \cdot \nabla \varphi - m \int_{M} \varphi^\theta v^{1-q}P^{m-1}\nabla v \cdot \nabla P  \,,
\end{multline*}
which, combined with \eqref{P_def12} and Cauchy Schwartz inequality, implies
\begin{equation} \label{vac2} 
c(n,q) \int_{M} v^{1-q}P^{m+1}\psi \leq  \theta \int_{M} \varphi^{\theta-1} v^{1-q}P^m |\nabla v||\nabla \varphi| + |m| \int_{M} \varphi^\theta v^{1-q}P^{m-1} |\nabla v| | \nabla P|  \,,
\end{equation}
where $c(n,q) >0$ provided that $q< \frac{n}{2}+1$. By Young inequality with exponents $(\gamma, \frac{\gamma}{\gamma-1})$ we bound the first term on the RHS of \eqref{vac2} as follows
\begin{equation*}
\theta \int_{M} \varphi^{\theta-1} v^{1-q}P^m |\nabla v||\nabla \varphi|  \leq \epsilon \int_M \vp^\theta v^{1-q} P^{m+1} + C(n,\theta,\epsilon,\gamma) \int_M \vp^{\theta -\gamma} v^{1-q} P^{m+1-\gamma} |\nabla v |^\gamma |\nabla \vp|^\gamma 
\end{equation*}
and, by using $|\nabla v|^2 < \frac{2}{n} vP$, we have
\begin{equation} \label{vac3}
\theta \int_{M} \varphi^{\theta-1} v^{1-q}P^m |\nabla v||\nabla \varphi|  \leq \epsilon \int_M \vp^\theta v^{1-q} P^{m+1} + \frac{C(n,\theta,\epsilon,\gamma)}{R^\gamma} \int_{B_{2R} \setminus B_R} \vp^{\theta -\gamma} v^{1-q+ \frac{\gamma}{2}} P^{m+1-\frac{\gamma}{2}} \,,
\end{equation}
where we also used the definition of $\vp$ \eqref{vphi_def}. 

Now we consider the second term on the RHS of \eqref{vac2} by using Young's inequality:
\begin{equation*}
|m| \int_{M} \varphi^\theta v^{1-q}P^{m-1} |\nabla v| | \nabla P| \leq \epsilon R^2 \int_{M} \varphi^{\theta+2} P^{t-2} v^{2-n} | \nabla P|^2 + \frac{C(\epsilon, m)}{R^2} \int_{M} \varphi^{\theta-2} v^{n-2q}P^{2m-t} |\nabla v|^2 
\end{equation*}
and using again $|\nabla v|^2 < \frac{2}{n} vP$ we obtain
\begin{equation} \label{vac4}
|m| \int_{M} \varphi^\theta v^{1-q}P^{m-1} |\nabla v| | \nabla P| \leq \epsilon R^2 \int_{M} \varphi^{\theta+2} P^{t-2} v^{2-n} | \nabla P|^2 + \frac{C(\epsilon, m,n)}{R^2} \int_{M} \varphi^{\theta-2} v^{1+n-2q}P^{2m-t+1} \,.
\end{equation}
The conclusion immediately follows from \eqref{vac2}, \eqref{vac3} and \eqref{vac4} by choosing $\epsilon < \frac12 c(n,q)$, where $c(n,q)$ is the constant appearing on the LHS of \eqref{vac2}.
\end{proof}

\medskip 

Now we are ready to give the proof of Theorem \ref{thm_critical}.

\begin{proof}[Proof of Theorem \ref{thm_critical}] The main goal of the proof is to show that $P$ is constant, i.e. that $|\nabla P| \equiv 0$ in $M$. Indeed, if $P$ is constant then \eqref{eqfond} implies that $k[v] \equiv 0$ in $M$ and, since $\Delta v = P = const$, \eqref{k_def2} yields 
$$
\nabla^2 v \equiv \frac{P}{d} g \quad  \textmd{ and } \quad  \Ric(\nabla v, \nabla v) \equiv 0 
$$
in $M$. From a classical splitting theorem (see \cite[Theorem 5.7.4]{Petersen} and also \cite{FMM} and \cite{CatinoMonticelli}) we obtain that $(M,g)$ is isometric to a complete flat warped product metric and hence it must be the Euclidean space of dimension $d$. Moreover, we immediately obtain that $v$ is quadratic, which implies that $u$ is given by  \eqref{u_radial_crit}.

Hence, the rest of the proof is devoted to show that $\nabla P \equiv 0 $ in $M$. 

%\end{proof}

%\begin{proposition} \label{prop_A}
%Let $(M,g)$ be a Riemannian manifold of dimension $d\geq 3$ with $\Ric \geq 0$. Let $P$ and $n$ be given by %\eqref{P_def12} and \eqref{n_def}, respectively. Assume that 
%\begin{equation} \label{n_cond}
% 3 \leq n \leq 5 \,.
%\end{equation}
%Then $P$ is constant. 
%\end{proposition}

\bigskip

\noindent $(i)$ 
We first observe that $n=d \geq 3$ in this case and hence Proposition \ref{mainineq} applies, and we have that 
\beq\label{comp26_bis} 
\left(t-\frac12\right)^2 \int_M  \vp^\theta  P^{t-2}v^{2-n}|\nabla P|^2  \leq \frac{4}{R^2} \int_{B_{2R} \setminus B_R} \vp^{\theta -2 } v^{2-n}P^t \,,
\eeq
for any $t>\frac{1}{2}$ and $R>1$. 

Now, we divide the proof in three cases accordingly to the value of $n$. 

\medskip

\noindent \emph{-- Case 1: $n =3$}. In this case we just set $t=1$ in \eqref{comp26_bis} and Lemma \ref{Lemma2bis} (we apply part (i) with $q=n-1$, so that $2 \leq q < \frac{n}{2} + 1$) yields that 
$$
 \int_{B_{2R} } v^{2-n}P = \int_{B_{2R} } v^{-(n-1)}\left(\frac{n}{2} |\nabla v|^2 + c_n \right)  \leq 2^{d+1} \omega_d C R^{d-n+1} 
$$
and from \eqref{comp26_bis} we immediately obtain that 
$$
\int_{B_R}  P^{-1}v^{2-n}|\nabla P|^2  \leq \frac{2^{d+5} \omega_d C}{R} \,.
$$
By letting $R \to +\infty$ we find that $|\nabla P| \equiv 0$ in $M$.

\noindent \emph{-- Case 2: $n=4$}. Again, we have to estimate the RHS of \eqref{comp26_bis}. To this end we choose 
\begin{equation} \label{eq_123}
t=\frac12 + \delta \quad  \textmd{ with } \quad \delta = (5-n)/5 \in \left( 0,\frac{1}{5}\right) 
\end{equation}
and set
\begin{equation} \label{eq_1234}
\gamma=\frac{n}{2} - \delta \,.
\end{equation}
%Since $t>1$ we can use H\"older's inequality 
%in order to apply Lemma \ref{Lemma2bis}. By setting 
%\begin{equation} \label{eq_1234}
%\gamma=\frac{n}{2} - \delta \,.
%\end{equation}
Since $ \frac{1}{2}<t<1$ we can use H\"older's inequality with parameters $(\frac{1}{t},\frac{1}{1-t})$ to get 
\begin{multline*}
\int_{B_{2R}} v^{2-n} P^t  =  \int_{B_{2R}} v^{2-n + \gamma t} v^{-\gamma t} P^t  \leq  \left(\int_{B_{2R}} v^{\frac{2-n + \gamma t}{1-t}}  \right)^{1-t} \left( \int_{B_{2R}} v^{- \gamma }  P \right)^t \\
=  \left(\int_{B_{2R}} v^{-\frac{n-2 - \gamma t}{1-t}}  \right)^{1-t} \left( \int_{B_{2R}} v^{- (1+\gamma) }\left(\frac{n}{2} |\nabla v|^2 + c_n\right) \right)^t
\end{multline*}
and Lemma \ref{Lemma2bis} yields
\begin{equation} \label{eq_12}
\frac{1}{R^2} \int_{B_{2R}} v^{2-n} P^t  \leq C R^{d-n-t} \,,
\end{equation}
provided that
\begin{equation} \label{12345}
\begin{cases}
\frac12 < t < 1 & \\
2 \leq 1+ \gamma <  \frac{n}{2}+1 & \\
0 \leq \frac{n-2 - \gamma t}{1-t}  \leq \frac{n}{2}+1 \,. & 
\end{cases}
\end{equation}
We notice that the first two lines in \eqref{12345} are verified thanks to \eqref{eq_123} and \eqref{eq_1234}. A simple calculation shows that also the last line in \eqref{12345} is satisfied. Hence from \eqref{eq_12} and \eqref{comp26_bis} we obtain that $|\nabla P| \equiv 0 $.

\noindent \emph{-- Case 3: $n=5$}. Let $\delta$ be a small parameter to be chosen later. Since $n=5$ then $\frac{n}{2} - 2 + \delta >0 $ and from $1<vP$ we can use 
$$
1 < (vP)^{\frac{n}{2} - 2 + \delta}
$$
in the integral on the RHS of \eqref{comp26_bis}, and obtain that 
\begin{equation} \label{n5_1}
\frac{1}{R^2} \int_{B_{2R} \setminus B_R} \vp^{\theta -2 }v^{2-n}P^t \leq \frac{1}{R^2} \int_{B_{2R} \setminus B_R} \vp^{\theta -2 } v^{-\frac{n}{2} +\delta}P^{t+ \frac{n}{2} - 2 + \delta} \,. 
\end{equation}
We apply Proposition \ref{prop_ngrande}, by letting 
\begin{equation} \label{param1}
m+1=t+ \frac{n}{2} -2 +\delta \quad \textmd{ and } \quad q=\frac{n}{2} + 1 - \delta
\end{equation}
(and $\theta-2$ in place of $\theta$) and the other parameters to be chosen later. From \eqref{comp26_bis}, \eqref{n5_1} and Proposition \ref{prop_ngrande} we obtain that 
\begin{multline*} 
\left(\Big(t-\frac12\Big)^2 - \epsilon\right) \int_M  \vp^\theta  P^{t-2}v^{2-n}|\nabla P|^2  \\
\leq \frac{C}{R^{2+\gamma}} \int_{B_{2R} \setminus B_R} \varphi^{\theta-2} v^{-\frac{n}{2} + \delta + \frac{\gamma}{2}} P^{t + \frac{n-4}{2} + \delta - \frac{\gamma}{2}} + \frac{C}{R^4} \int_{B_{2R}} \varphi^{\theta-4} v^{-1 + 2 \delta} P^{t+n-5+2 \delta} \,,
\end{multline*}
and by H\"older's inequality with exponents $(s,\frac{s}{s-1})$ and since $\supp \vp \subset B_{2R}$, we obtain
\begin{multline}\label{n5_10}
\left(\Big(t-\frac12\Big)^2 - \epsilon\right) \int_M  \vp^\theta  P^{t-2}v^{2-n}|\nabla P|^2  \leq \frac{C}{R^{2+\gamma}} \int_{B_{2R} \setminus B_R} v^{-\frac{n}{2} + \delta + \frac{\gamma}{2}} P^{t + \frac{n-4}{2} + \delta - \frac{\gamma}{2}} \\ 
+\frac{C}{R^4} \left( \int_{B_{2R}} v^{-s + 3s\delta} P^{(t  + n-5 + 2\delta)s} \right)^{\frac{1}{s}} \left( \int_{B_{2R}} v^{- \frac{s}{s-1}\delta} \right)^{\frac{s-1}{s}} 
 \,.
\end{multline}
Now we have to choose the parameters in order to apply Lemma \ref{Lemma2bis} (see also Remark \ref{remark_asymp}). We choose $t> \frac12$, $s>1$ and $\delta>0$ small such that
\begin{equation} \label{para_1}
\begin{cases}
-\frac{n}{2} < -s + 3s\delta \leq -1 & \\
(t  + n-5 + 2\delta)s = 1 \,.& 
\end{cases}
\end{equation}
For instance, since $n=5$, by choosing $t = \frac12 + \delta$ and $s= (\frac12 + 2 \delta)^{-1}$ we have that the conditions in \eqref{para_1} are satisfied. Since $s,\gamma >1$ then 
$$
t + \frac{n-4}{2} + \delta - \frac{\gamma}{2} < (t + n-5+ 2\delta )s = 1
$$
which implies 
$$
1- t -  \frac{n-4}{2} - \delta + \frac{\gamma}{2} >0 \,,
$$
and from $1 < vP$ we have that
$$
1 \leq (vP)^{1-t - \frac{n-4}{2} - \delta + \frac{\gamma}{2}}
$$
which we use in the first integral on the RHS in \eqref{n5_10}. Hence, from \eqref{n5_10} we obtain 
\begin{multline}\label{n5_101}
\left(\Big(t-\frac12\Big)^2 - \epsilon\right) \int_M  \vp^\theta  P^{t-2}v^{2-n}|\nabla P|^2  \leq \frac{C}{R^{2+\gamma}} \int_{B_{2R} \setminus B_R} v^{-n + 3 + \gamma - t } P  \\ 
+\frac{C}{R^4} \left( \int_{B_{2R}} v^{-s + 3s\delta} P  \right)^{\frac{1}{s}} \left( \int_{B_{2R}} v^{- \frac{s}{s-1}\delta} \right)^{\frac{s-1}{s}} 
 \,.
\end{multline}
By choosing $\gamma>1$ such that $ -\frac{n}{2} < -n + 3 + \gamma - t \leq -1$, i.e., 
$$
\frac{n}{2} - 3 + t < \gamma \leq n-4 +t \,,
$$
$\delta>0$ small enough and from \eqref{para_1}, we apply Lemma \ref{Lemma2bis} and from \eqref{n5_101} we obtain 
$$
\left(\Big(t-\frac12\Big)^2 - \epsilon\right) \int_M  \vp^\theta  P^{t-2}v^{2-n}|\nabla P|^2  \leq C R^{d -n - t} + C R^{d -5  - \frac{1}{s} + 2 \delta} \,.
$$
Once $t> \frac12$ is fixed, since $n=d=5$ we obtain that 
$$
\lim_{R\to \infty} \int_{B_R} P^{t-2}v^{2-n}|\nabla P|^2 = 0 \,,
$$
which implies that $P$ is constant in $M$.

\bigskip

Now we turn to the proof of item $(iii)$. 

\bigskip

\noindent $(iii)$ We first prove that the assumption $u \in L^{\frac{2d}{d-2}}(M)$ implies $\int_M \vert \nabla u\vert^2 < + \infty $. To this end, we multiply the equation \eqref{eq_crit} by $u \varphi^2$, where $\varphi$ is a standard cutoff function as in \eqref{vphi_def}, and then we integrate by parts to get 
\begin{equation}
\int_M \vert \nabla u\vert^2 \varphi^2 = \int_M u^{\frac{2d}{d-2}} \varphi^2 - \int_M 2 u \varphi \nabla u \cdot \nabla \varphi 
\leq \int_M u^{\frac{2d}{d-2}} \varphi^2 + \frac{1}{2}\int_M \vert \nabla u\vert^2 \varphi^2 + 2 \int_M  u^2 \vert  \nabla \varphi \vert^2.  
\end{equation}
Hence, for every $R>0$, 
\begin{equation*}
\begin{split}
\frac{1}{2} \int_{B_R} \vert \nabla u\vert^2 & \leq  \frac{1}{2} \int_M \vert \nabla u\vert^2 \varphi^2 \leq \int_M u^{\frac{2d}{d-2}} \varphi^2 + 2 \int_M  u^2 \vert  \nabla \varphi \vert^2 \\
& \leq \int_M u^{\frac{2d}{d-2}} \varphi^2 + 2 \left( \int_M u^{\frac{2d}{d-2}} \right)^{\frac{d-2}{d}}  
\left( \int_M  \vert \nabla \varphi \vert^d \right)^{\frac{2}{d}} \\
& \leq   \int_M u^{\frac{2d}{d-2}}  + 2 \left( \int_M u^{\frac{2d}{d-2}} \right)^{\frac{d-2}{d}} \left( R^{-d} V_g(B_{2R}) \right)^{\frac{2}{d}} 
\leq  \int_M u^{\frac{2d}{d-2}}  + 8 \omega_d^{\frac{2}{d}}\left( \int_M u^{\frac{2d}{d-2}} \right)^{\frac{d-2}{d}}
\end{split}
\end{equation*}
where in the latter we have used  \eqref{BG}. Then, we obtain $\int_M \vert \nabla u\vert^2 < + \infty $ by letting $ R \longrightarrow \infty$.

Next, we fix a point $p \in M$ and apply Lemma \ref{lemma_Pt} with $t=1$, $\theta \geq 2$ and  $\psi$ a standard cutoff function as in \eqref{vphi_def}. This leads to 

\begin{multline} \label{stima-ast0}
\frac{1}{2}  \int_{M}  P^{-1}v^{2-d}|\nabla P|^2\psi^\theta  \le - \theta\int_{M} \psi^{\theta-1}  v^{2-d} \nabla P \cdot \nabla\psi\\
\le \theta \left(\int_{M} \psi^\theta v^{2-d}P^{-1}|\nabla P|^2\right)^\frac{1}{2}\left(\int_{M} \psi^{\theta-2}v^{2-d}P \vert \nabla \psi|^2\right)^\frac{1}{2}.
\end{multline}

Then, we observe that 

\begin{multline} \label{stima-ast}
\int_{M} \psi^{\theta-2}v^{2-d}P \vert \nabla \psi|^2 \leq \frac{1}{R^2} \int_{B_{2R} \setminus B_R} v^{2-d} P = \frac{1}{R^2} \int_{B_{2R} \setminus B_R}  v^{1-d} Pv \\
= \frac{1}{R^2} \int_{B_{2R} \setminus B_R}  v^{1-d} \left( \frac{d}{2} \vert \nabla v|^2 + \frac{2}{d-2}\right)  =
\frac{1}{R^2} \int_{B_{2R} \setminus B_R}  v \left( \frac{d}{2} v^{-d} \vert \nabla v|^2 + \frac{2}{d-2} v^{-d} \right) \\
= \frac{1}{R^2} \int_{B_{2R} \setminus B_R}  v \left( \frac{2d}{(d-2)^2}  \vert \nabla u|^2 + \frac{2}{d-2}  u^{\frac{2d}{d-2}} \right)
\end{multline}
and we recall that, since $u$ is superharmonic on $M$ and $\Ric \geq 0$, we have the existence of two positive constants $A_0, R_0$ such that 
\begin{equation} 
u(x) \geq \frac{A_0}{r(x)^{d-2}}\qquad \forall \, x \in M \setminus B_{R_0},
\end{equation}
where $r=r(x)$ denotes the geodesic distance from a fixed reference point of $M$ 
(see for instance \cite{FMM} and \cite{CatinoMonticelli}). The latter implies that 
\begin{equation} 
v(x) \leq A_0^{\frac{2}{2-d}} r(x)^{2} \qquad \forall \, x \in M \setminus B_{R_0}
\end{equation}
and therefore, from \eqref{stima-ast}, we get that, for every $R> R_0$  

\begin{multline}\label{ultimaa}
\int_{M} \psi^{\theta-2}v^{2-d}P \vert \nabla \psi|^2 \leq 
\frac{1}{R^2} 
\int_{B_{2R} \setminus B_R}  v \left( \frac{2d}{(d-2)^2}  \vert \nabla u|^2 + \frac{2}{d-2}  u^{\frac{2d}{d-2}} \right)
\\ \leq 2d A_0^{\frac{2}{2-d}} \int_{B_{2R} \setminus B_R} \left( \vert \nabla u|^2 + u^{\frac{2d}{d-2}} \right).
\end{multline}

Now we plug \eqref{ultimaa} into \eqref{stima-ast0} to get 

\begin{multline} 
\int_{B_R}  P^{-1}v^{2-d}|\nabla P|^2 \leq \int_{M}  P^{-1}v^{2-d}|\nabla P|^2\psi^\theta  \leq 4 \theta^2 \int_{M} \psi^{\theta-2}v^{2-d}P \vert \nabla \psi|^2 \\
\leq 8 d \theta^2 A_0^{\frac{2}{2-d}} \int_{B_{2R} \setminus B_R} \vert \nabla u|^2 + u^{\frac{2d}{d-2}}\longrightarrow 0 \qquad as \quad R \longrightarrow \infty,
\end{multline}
since $u \in L^{\frac{2d}{d-2}}(M)$ and $\int_M \vert \nabla u\vert^2 < + \infty $.

\bigskip

Finally, we prove item $(ii)$. 

\bigskip

\noindent $(ii)$ If $\alpha \geq \frac{d-2}{2}$, then $u \in L^{\frac{2d}{d-2}}(M)$, and so the result follows 
from item $(iii)$. We can therefore assume $\alpha < \frac{d-2}{2}$.

\noindent Our goal is to estimate the RHS in \eqref{comp26_bis}. If we assume that $t \geq 1$, then the RHS in \eqref{comp26_bis} can be estimated by
\begin{equation} \label{iii-1} 
\frac{4}{R^2} \int_{B_{2R} \setminus B_R} \vp^{\theta -2 } v^{2-d}P^t \leq \frac{4 }{R^2} 2^{t-1}\left( \frac{d^t}{2^t} \int_{B_{2R} \setminus B_R} v^{2-d-t} |\nabla v|^{2t} +  \frac{2^t}{(d-2)^t}   \int_{B_{2R} \setminus B_R} v^{2-d-t}  \right) \,.
\end{equation}

Let
$$
t= \frac{d-2}{2} 
$$
and let $\delta>0 $ to be chosen sufficiently small. 

We consider the first integral on the RHS of \eqref{iii-1}, and we write
\begin{equation} \label{eq_gian}
v^{2-d-t} |\nabla v|^{2t} = v^{-\frac{3}{2}(d-2)} |\nabla v|^{d-2} =     v^{-\delta} \left(\frac{|\nabla v|}{v} \right)^{d-4}  
v^{-(\frac{d}{2}+1-\delta)} |\nabla v|^2\,.
\end{equation}
Recalling that $u=v^{-\frac{d-2}{2}}$, from \cite[Proposition 2.1]{FMM}, we have that 
\begin{equation} \label{eq_FMM}
\sup_{B_\rho(x)} \frac{|\nabla v|^2}{v^2} \leq C \left(\frac{1}{\rho^2} + \sup_{B_{2\rho} (x)} \frac{1}{v^{2}} \right)
\end{equation}
for any $x \in M$, $\rho>1$ and for some dimensional constant $C>0$.
Since $u(x) = O(r^{-\alpha})$ as $r \to \infty$, then 
\begin{equation} \label{v-1}
v^{-1} = O(r^{-\frac{2\alpha}{d-2}}) \,,
\end{equation} 
and from \eqref{eq_FMM} we obtain that\footnote{This can be obtained by choosing $\rho = R/4$ in \eqref{eq_FMM} and arguing as in the proof of \cite[Corollary 2.2]{FMM}.}
\begin{equation} \label{eq_yau}
\frac{|\nabla v|^2}{v^2} \leq C R^{-\frac{4\alpha}{d-2}}
\end{equation}
in $B_{2R} \setminus B_R$, as $R \to \infty$. Hence, from \eqref{eq_gian}, \eqref{v-1} and \eqref{eq_yau} we find
\begin{equation*}
\frac{1}{R^2} \int_{B_{2R} \setminus B_R} v^{2-d-t} |\nabla v|^{2t}  \leq C R^{-2-\frac{2\alpha \delta}{d-2} -\frac{2\alpha}{d-2}(d-4) }  \int_{B_{2R} \setminus B_R} v^{- (\frac{d}{2}+1-\delta)} |\nabla v|^{2} \,.
\end{equation*}
By applying Lemma \ref{Lemma2bis}, with $q= \frac{d}{2}+1-\delta$ and any $ \delta \in \left( 0, \frac{d-2}{2} \right]  $, we obtain
\begin{equation} \label{prim_int}
\frac{1}{R^2} \int_{B_{2R} \setminus B_R} v^{2-d-t} |\nabla v|^{2t}  \leq C R^{-\frac{2\alpha }{d-2} (d-4 + \delta) + \frac{d-6}{2} + \delta }  \,,
\end{equation}
as $R\to \infty$.

We now consider the second integral on the RHS of \eqref{iii-1}. Since 
$$
v^{2-d-t} = v^{-\frac32(d-2)} = v^{-(d-4)} v^{-(\frac{d}{2} + 1)} 
$$ 
and by using \eqref{v-1}, the second integral on the RHS of \eqref{iii-1} can be estimated by 
$$
\frac{1}{R^2} \int_{B_{2R} \setminus B_R} v^{2-d-t} \leq C R^{-2 - \frac{2\alpha}{d-2}(d-4)}  \int_{B_{2R} \setminus B_R} 
v^{-(\frac{d}{2} + 1)} 
$$
as $R \to \infty$, and Lemma \ref{Lemma2bis}, with $q= \frac{d}{2}+1$, yields
\begin{equation} \label{sec_int}
\frac{1}{R^2} \int_{B_{2R} \setminus B_R} v^{2-d-t} \leq  C R^{- \frac{2\alpha}{d-2}(d-4) + \frac{d-6}{2}} \,.
\end{equation}
From \eqref{iii-1}, \eqref{prim_int}, \eqref{sec_int} and the assumption $\alpha > \frac{(d-2)(d-6) }{4(d-4)}$, we can choose a positive and small enough $\delta$ such that
$$
\lim_{R\to \infty} \frac{4}{R^2} \int_{B_{2R} \setminus B_R} \vp^{\theta -2 } v^{2-d}P^t = 0 
$$
and from \eqref{comp26_bis} we obtain that $P$ is constant, which completes the proof.

\end{proof}

\section{Classification results for Liouville equation} \label{Section_Liouville}
In this section we prove Theorem \ref{thm_Liouville} and Theorem \ref{thm_Liouville2} where solutions to the Liouville equation 
$$
-\Delta u = e^u  \qquad \textmd{ in } \quad M 
$$
are considered. In this case, we recall that $n=d=2$ and that, by letting 
$$
u=-2 \log v \,,
$$ 
we find that $v$ satisfies
$$
\Delta v = P \,,
$$
where
$$
P = \frac{1}{v} \left( |\nabla v|^2 + \frac12 \right)  \,.
$$

\begin{proof}[Proof of Theorem \ref{thm_Liouville}] 
By applying Lemma \ref{Lemma1} with $ q=1$ and $n=d=2$ we get 
\begin{equation} \label{eq_Lemma1-exp}
%\int_{M} P \psi = 
\int_{M} v^{-1}|\nabla v|^2\psi+ \frac{1}{2}\int_{M} v^{-1}\psi = 
- \int_{M} \nabla v \cdot \nabla \psi
\end{equation}
for any $\psi\in  C^{0,1}_c(M)$. By using $\psi = \varphi^2$ in \eqref{eq_Lemma1-exp}, where $\varphi$ is a standard cutoff function as in \eqref{vphi_def}, we deduce 

\begin{multline} \label{stima1-exp}
\int_{M} v^{-1}|\nabla v|^2\varphi^2+ \frac{1}{2}\int_{M} v^{-1}\varphi^2 = - \int_{M} \nabla v \cdot \nabla \varphi^2 = 
-2  \int_{M} \varphi  \nabla v \cdot \nabla \varphi \\
\leq \int_{M} \left( \varphi v^{-\frac{1}{2}} \vert \nabla v \vert \right) \left(2 v^{\frac{1}{2}} \vert \nabla \varphi \vert \right)
\leq \frac{1}{2}\int_{M} v^{-1}|\nabla v|^2\varphi^2 + 2 \int_{M} v \vert \nabla \varphi \vert ^2 \,,
\end{multline}
and so
\begin{equation} \label{stima2-exp}
\frac{1}{2}\int_{M}  P \varphi^2 \leq \frac{1}{2}\int_{M} v^{-1}|\nabla v|^2\varphi^2+ \frac{1}{2}\int_{M} v^{-1}\varphi^2 \leq 2 \int_{M} v \vert \nabla \varphi \vert ^2 \,.
\end{equation}
From assumption \eqref{eq_crit2-hyp} we infer that 
\begin{equation}
v(x) = e^{-\frac{u(x)}{2}}\leq r^2(x)G(r(x)), \qquad \textmd{for} \quad  r(x) >c,
\end{equation}
and therefore, for any $ R >c$, 
\begin{multline} \label{stima3-exp}
\int_{B_R} P  \leq 4 \int_{M} v \vert \nabla \varphi \vert ^2 = 4  \int_{B_{2R} \setminus B_R} v \vert \nabla \varphi \vert ^2 \leq 
\frac{4}{R^2} \int_{B_{2R} \setminus B_R} v \\
\leq  \frac{4}{R^2} (2R)^2 G(2R) V_g(B_{2R} \setminus B_R) \leq 4^3 \omega_2 R^2 G(2R), 
\end{multline}
where in the latter we have used \eqref{BG}, since $\Ric \geq 0$.

For $t >c$, set $ F(t) = G(2t)$, then  $ \int_c^{\infty} \frac{dt}{tF(t)} = \int_c^{\infty} \frac{2dt}{2tG(2t)}  = 
\int_{2c}^{\infty} \frac{dt}{tG(t)} = + \infty$, thanks to the integral condition satisfied by $G$. Hence, \eqref{stima3-exp} yields
\begin{equation} \label{stima4-exp}
\limsup_{R \to \infty} \frac{1}{R^2 F(R)} \int_{B_R} P \leq 4^3 \omega_2 .
\end{equation}

Now, by applying Lemma \ref{lemma_Pt} we get  $\Delta P^t \geq 0$ on $M$,  for any 
$ t \in \left[ \frac{1}{2}, 1\right).$ Since $ P^t \geq 0$ on $M$, we can thus apply \cite[Theorem 2.2]{Karp} to get that 
$P^t= const.$ on $M$ for some/any $t  \in \left[ \frac{1}{2}, 1\right).$ The latter implies that $P$ is constant on $M$ and the  desired conclusion then follows as in Theorem  \ref{thm_critical}.

\smallskip

The optimality of the coefficient $4$ in the lower bound \eqref{eq_crit2-hyp} can be deduced directly from the explicit examples provided on pp.12-13 of \cite{CaiLai}.
\end{proof}

\begin{proof}[Proof of Theorem \ref{thm_Liouville2}] 
Since $(M,g)$ is not compact, by the classical results of Cohn-Vossen \cite{CV}  and Huber \cite{HU} 
(see also \cite[Proposition 1.1]{LiTam}) we know that either $(M,g)$ is isometric to the flat cylinder $\mathbb{S}^1 \times \R$ or $(M,g)$  is conformal to $\R^2$ with the Euclidean metric. 
Therefore, to obtain the desired conclusion it is enough to prove that the assumption \eqref{eq_crit3-hyp} rules out the flat cylinder. 
For contradiction, suppose that $(M,g)$ is the flat cylinder, then $(M,g)$ has linear volume growth, i.e., there exists some  constant $C>0$ such that $V_g(B_R) \leq C R,$ for large radii $R$. Therefore, by proceeding as in the previous proof we have 

\begin{equation}
v(x) = e^{-\frac{u(x)}{2}}\leq r^3(x)G(r(x)), \qquad \textmd{for} \quad  r(x) >c,
\end{equation}
and so, for any $ R >c$, 
\begin{multline} \label{stima3-exp2}
\int_{B_R} P  \leq 4 \int_{M} v \vert \nabla \varphi \vert ^2 = 4  \int_{B_{2R} \setminus B_R} v \vert \nabla \varphi \vert ^2 \leq 
\frac{4}{R^2} \int_{B_{2R} \setminus B_R} v \\
\leq  \frac{4}{R^2} (2R)^3 G(2R) V_g(B_{2R} \setminus B_R) \leq 64 C R^2 G(2R), 
\end{multline}
where in the latter we have used the linear volume growth of the flat cylinder. 

As in the previous proof, inequality \eqref{stima3-exp2} allows us to conclude that $P$ is constant on $M$, and so $(M,g)$ should be isometric to $\R^2$ with the Euclidean metric. A contradiction. 

\smallskip

The optimality of the the above result  can be deduced directly from the explicit examples provided on pp.12-13 of \cite{CaiLai}.

\end{proof}

\medskip

\section*{Acknowledgments} 
\noindent The authors have been partially supported by the ``Gruppo Nazionale per l'Analisi Matematica, la Probabilit\`a e le loro Applicazioni'' (GNAMPA) of the ``Istituto Nazionale di Alta Matematica'' (INdAM, Italy) and by the Research Project of the Italian Ministry of University and Research (MUR) Prin 2022 ``Partial differential equations and related geometric-functional inequalities'', grant number 20229M52AS\_004.

After we posted on ArXiv a first version of the present paper we had many stimulating discussions with G. Catino, A. Ferrero, M. Fogagnolo, G. Grillo, A. Malchiodi, L. Mazzieri, D. Monticelli and M. Muratori.  We thank all of them since those discussions helped us to improve the presentation of the paper.  
We also express our gratitude to X. Cai and M. Lai, who kindly drew our attention to their very recent manuscript \cite{CaiLai}.

\end{document}